\documentclass[leqno]{article}
\usepackage{verbatim, amsmath, amscd, amssymb}
\usepackage[all]{xy}
\newcommand{\cI}{{\mathcal I}}
\newcommand{\cA}{{\mathcal A}}
\newcommand{\cB}{{\mathcal B}}
\newcommand{\cC}{{\mathcal C}}

\newcommand{\cF}{{\mathcal F}}

\newcommand{\cL}{{\mathcal L}}

\newcommand{\cO}{{\mathcal O}}
\newcommand{\cP}{{\mathcal P}}

\newcommand{\frU}{\mathfrak U}

\newcommand{\rC}{\mathrm{C}}

\newcommand{\rB}{\mathrm{B}}

\newcommand{\rF}{\mathrm{F}}\newcommand{\rG}{\mathrm{G}}\newcommand{\rH}{\mathrm{H}}

\newcommand{\rN}{\mathrm{N}}

\newcommand{\bbF}{\mathbb F}

\newcommand{\bbP}{\mathbb P}
\newcommand{\bbN}{\mathbb N}

\newcommand{\bbZ}{\mathbb Z}
\newcommand{\bfR}{\mathbf R}

\newcommand{\ad}{\mathrm{ad}}\newcommand{\Ad}{\mathrm{Ad}}

\newcommand{\Dist}{\mathrm{Dist}}

\newcommand{\ev}{\mathrm{ev}}

\newcommand{\Fr}{\mathrm{Fr}}
\newcommand{\Frac}{\mathrm{Frac}}

\newcommand{\id}{\mathrm{id}}
\newcommand{\im}{\mathrm{im}}

\newcommand{\ind}{\mathrm{ind}}

\newcommand{\Lie}{\mathrm{Lie}}

\newcommand{\Proj}{\mathrm{Proj}}

\newcommand{\res}{\mathrm{res}}

\newcommand{\SL}{\mathrm{SL}}

\newcommand{\St}{\mathrm{St}}

\newcommand{\Ab}{\mathbf{Ab}}

\newcommand{\Mod}{\mathbf{Mod}}
\newcommand{\Sch}{\mathbf{Sch}}

\newcommand{\lbr}{\begin{bmatrix}}
\newcommand{\rbr}{\end{bmatrix}}
\newcommand{\for}{\bigcirc\kern-2.6ex \because}
\newcommand{\forb}{\bigcirc\kern-2.8ex \because}
\newcommand{\forbb}{\bigcirc\kern-3.0ex \because}
\newcommand{\forbbb}{\bigcirc\kern-3.1ex \because}

\newcommand{\End}{\operatorname{End}}
\newcommand{\Ru}{\mathrm{Ru}}

\newcommand\pf{\noindent {\it D\'emonstration.  }}
\newtheorem{thm}{Th\'eor\`eme.}

\newtheorem{prop}{Proposition.}

\newtheorem{lem}{Lemme.}

\date{}
\begin{document}
\large
\title{
{\bf
Contraction par Frobenius de 
$G$-modules
}
\thanks
{supported in part
by JSPS Grant in Aid
for Scientific Research}
\author{
Michel G\textsc{ros}
\\
Universit\'e de Rennes I
\\
IRMAR
\\
michel.gros@univ-rennes1.fr
\and
K\textsc{aneda} Masaharu
\\
Osaka City University
\\
Department of Mathematics
\\
kaneda@sci.osaka-cu.ac.jp
}
}
\maketitle

\begin{abstract}
Let $G$ be a simply connected semisimple
algebraic group
over an algebraically closed field
$\Bbbk$
of positive characteristic.
We will untwist the structure of 
$G$-modules by a newly found splitting of the Frobenius endomorphism on the algebra of distributions of $G$,
which, for example, explains the subsumed $G$-equivariant nature of 
the Frobenius splitting of the flag variety of $G$.

\end{abstract}

\begin{center}
{\bf
Introduction}
\end{center}

Soient $\Bbbk$ un corps alg\'ebriquement clos  de caract\'eristique $p>0$ et
$G_\Bbbk$ un groupe alg\'ebrique semi-simple simplement connexe.
Soit $M$ un $G_\Bbbk$-module ;
l'action de $G_\Bbbk$ via le Frobenius fournit sur  $M$ une nouvelle structure de $G_\Bbbk$-module que nous noterons  $M^{[1]}$.
Si, de plus, le noyau de Frobenius  de   $G_\Bbbk$ agit trivialement sur $M$,
$M$ admet alors une autre structure de   $G_\Bbbk$-module not\'ee
$M^{[-1]}$ telle que si on lui applique  le twist de Frobenius pr\'ec\'edent  
$(M^{[-1]})^{[1]}$, on retrouve  $M$.
En utilisant une g\'en\'eralisation du scindage de Frobenius sur  l'alg\`ebre des distributions  $\Dist(G_\Bbbk)$ de $G_\Bbbk$
d\'ecrit  dans \cite{G}, on peut "d\'etwister" l'action de Frobenius de n'importe quel    $G_\Bbbk$-module.
Plus pr\'ecis\'ement, il existe une mesure invariante involutive   $\mu_0$ dans l'alg\`ebre des distributions du noyau de Frobenius d'un tore maximal $T_\Bbbk$ de $G_\Bbbk$ et un homomorphisme d'alg\`ebres $\phi$ de $\Dist(G_\Bbbk)$ vers une "sous-alg\`ebre" de $\Dist(G_\Bbbk)$ avec  $\mu_0$ comme unit\'e de ce dernier qui scinde l'endomorphisme de Frobenius sur $\Dist(G_\Bbbk)$.
Si $\Lambda$ d\'esigne le groupe des poids de
$T$
et $M_\lambda$,
$\lambda\in\Lambda$,
le sous-espace de poids  $\lambda$ de
$M$, alors
$\mu_0M=\coprod_{\lambda\in\Lambda}M_{p\lambda}$,
et cela munit donc ce dernier d'une  $G_\Bbbk$-action via $\phi$,
que nous noterons
$M^\phi$,
telle que l'on ait 
$(M^{[1]})^\phi\simeq M$. 
Ainsi  $M^\phi$ est la contraction de $M$ que d\'efinissait
Littelmann  \cite{Li}
en utilisant  le morphisme de Frobenius quantique introduit par  Lusztig  \cite{L90}.

D'autre part, Kumar et Littelmann \cite{KuLi} ont \'etendu le scindage de ce morphisme, d\'efini initialement seulement sur la partie unipotente  \cite{L90} \`a la partie Borel et ont appliqu\'e   ceci \cite{KuLi} \`a la cohomologie des fibr\'es en droites sur la vari\'et\'e  des drapeaux   $G_\Bbbk/B_\Bbbk$,
(avec $B_\Bbbk$ un sous-groupe de  Borel de $G_\Bbbk$).
La sp\'ecialisation en 1 du param\`etre quantique $q$ leur  donne alors une construction alg\'ebrique du scindage de Frobenius de   $G_\Bbbk/B_\Bbbk$
\cite{MR}, qui est  $B_\Bbbk$-\lq\lq
semi-invariant" (nous adoptons ici la terminologie sugg\'er\'ee par
\cite[4.1.5]{BK},  plut\^ot que la terminologie ``canonique").

Le scindage 
$\phi$ sur 
$\Dist(G_\Bbbk)$ que nous \'etudierons permet de  reformuler \cite{KuLi} de mani\`ere   $G$-\'equivariante et par suite de rendre   plus naturelle la  propri\'et\'e de $B_\Bbbk$-semi-invariance du scindage de Frobenius de  $G_\Bbbk/B_\Bbbk$. Le scindage obtenu est \'egalement  $B^+_\Bbbk$-semi-invariant (mais ne peut \^etre rendu  $G_\Bbbk$-\'equivariant).
Nos arguments sont une variante modulaire, souvent plus simple,  de ceux de 
\cite{KuLi} et certaines preuves formellement analogues \`a celles de loc. cit. seront souvent omises. L'existence d'un scindage de Frobenius sur un sch\'ema ayant des corollaires importants (d\'ej\`a largement d\'egag\'es par Mehta,  Ramanathan, Ramanan, Mathieu,...), nous esp\'erons que cette reformulation m\'eritait d'\^etre faite.

Dans le premier paragraphe, nous d\'efinissons le scindage de Frobenius sous des hypoth\`{e}ses tr\`es g\'en\'erales (le cas des ``petites"caract\'eristiques n\'ecessitant un peu de travail), puis dans les deux suivants, comment ce dernier permet de r\'einterpr\'eter en toute caract\'eristique les constructions de \cite{Li} et de \cite{KuLi}. Enfin, dans le dernier paragraphe, nous expliquons comment faisceautiser ces constructions et comment elles permettent \'egalement de v\'erifier que les sous-sch\'emas de Schubert sont  scind\'es de mani\`ere compatible.

Le pr\'esent travail s'est d\'evelopp\'e pendant la visite du second auteur \`a l'IRMAR au printemps 2009. Il remercie  cet institut pour son hospitalit\'e et son support financier.

\begin{center}
{\bf
$1^\circ$
Scindage de Frobenius sur $\Dist(G)$}
\end{center}

Pour simplifier les notations, nous travaillerons, sauf mention du contraire, sur le corps $\bbF_p$ \`a  $p$ \'el\'ements.
Soit donc $G$ un groupe alg\'ebrique semi-simple simplement connexe sur
$\bbF_p$,
$B$ un sous-groupe de Borel de $G$, et $T$
un tore maximal de  $B$ tous deux scind\'es sur $\bbF_p$.
Soit
$\Dist(G)$ (resp. $\Dist(B)$,
$\Dist(T)$)
l'alg\`ebre des distributions de
$G$
(resp. $B$, $T$).
Soit
$\Lambda$ le groupe des poids de
$B$, $R\subset\Lambda$ l'ensemble des racines de $G$ relativement \`a $T$
et
$R^+$ un syst\`eme positif
tel que les racines de  $B$ soient n\'egatives ; $R^s=\{\alpha_i\mid
i\in[1,\ell]\}$ l'ensemble des racines simples.
Pour  $\alpha\in R$, nous noterons  $\alpha^\vee$
la co-racine de
$\alpha$.
Si $M$ est un $T$-module
et si $\lambda\in\Lambda$,
$M_\lambda$
d\'enotera le sous-espace de poids  $\lambda$ de $M$.

Soit
$Fr : G\to
G$
l'endomorphisme de Frobenius de $G$.
Le  comorphisme correspondant  $Fr^\sharp : \bbF_p[G]\to\bbF_p[G]$, donn\'e par 
 $a\mapsto a^p$,
induit un homomorphisme de
$\bbF_p$-alg\`ebres  
$\Dist(Fr) : \Dist(G)\to\Dist(G)$
via
$\mu\mapsto \mu\circ
Fr^\sharp$.

Nous allons \'etendre au cas g\'en\'eral la construction du scindage de Frobenius sur   $\Dist(G)$
consid\'er\'ee dans  \cite{G} pour $G=\SL_2$.

\noindent
(1.1) Rappelons tout d'abord quelques r\'esultats sur les alg\`ebres de distributions. 
Soit
$G_\bbZ$
la
$\bbZ$-forme de Chevalley  de
$G$. Soit $E_i$ (resp. $F_i$), 
$i\in[1,\ell]$,
l'\'el\'ement d'une base de Chevalley de
$\Lie(G_\bbZ)$
correspondant \`a la racine simple
$\alpha_i$ (resp.
$-\alpha_i$). Posons
$H_i=[E_i,F_i]$. Soit aussi $U_\bbZ$
(resp. $U_\bbZ^+$)
la $\bbZ$-forme de
$U=\Ru(B)$
(resp. son oppos\'ee
$U^+$).
On a alors  (\cite[Satz I.7]{J73})
\begin{align}
\Dist(U_\bbZ)
&=
\bbZ[F_i^{(k)}\mid
i\in I, k\in\bbN],
\quad
\Dist(U_\bbZ^+)=
\bbZ[E_i^{(k)}\mid
i\in I, k\in\bbN],
\\
\notag
\Dist(G_\bbZ)
&=\bbZ[
E_i^{(n)},F_i^{(n)}
\mid
i\in[1,\ell],k\in\bbN].
\end{align}

Pour v\'erifier ces \'egalit\'es,   signalons qu'on peut aussi utiliser l'action du groupe de Weyl sur 
  $\Dist(G_\bbZ)$. Si
$s_i=x_{\alpha_i}(1)x_{-\alpha_i}(-1)x_{\alpha_i}(1)$ (notations comme dans  \cite[II.1.2]{J}), alors via la repr\'esentation adjointe

 \[
s_iE_j=
\begin{cases}
E_j
&\text{si $\langle\alpha_j,\alpha_i^\vee\rangle=0$},
\\
-F_i
&\text{si $j=i$
(et alors  $\langle\alpha_j,\alpha_i^\vee\rangle=2$)},
\\
E_iE_j-E_jE_i
&\text{si $\langle\alpha_j,\alpha_i^\vee\rangle=-1$},
\\
E_i^{(2)}E_j-E_iE_jE_i+E_jE_i^{(2)}
&\text{si $\langle\alpha_j,\alpha_i^\vee\rangle=-2$},
\\
E_i^{(3)}E_j-E_i^{(2)}E_jE_i+E_iE_jE_i^{(2)}-E_jE_i^{(3)}
&\text{si $\langle\alpha_j,\alpha_i^\vee\rangle=-3$},
\end{cases}
\]
et donc aussi, en utilisant 
l'anti-automorphisme $\Omega$  de $\Dist(G_\bbZ)$
introduit dans  [L90, 1.1 (d1)],
\[
s_iF_j=
\begin{cases}
F_j
&\text{si $\langle\alpha_j,\alpha_i^\vee\rangle=0$},
\\
-E_i
&\text{si $j=i$
(et alors   $\langle\alpha_j,\alpha_i^\vee\rangle=2$)},
\\
F_jF_i-F_iF_j
&\text{si $\langle\alpha_j,\alpha_i^\vee\rangle=-1$},
\\
F_jF_i^{(2)}-F_iF_jF_i+F_i^{(2)}F_j
&\text{si $\langle\alpha_j,\alpha_i^\vee\rangle=-2$},
\\
F_jF_i^{(3)}+F_i^{(2)}F_jF_i-F_iF_jF_i^{(2)}-F_i^{(3)}F_j
&\text{si $\langle\alpha_j,\alpha_i^\vee\rangle=-3$}.
\end{cases}
\]
 On en d\'eduit alors comme dans 
\cite[41.1.3]{L93}/\cite[5.7]{L90} 
que
$
\Dist(U_\bbZ)=
\bbZ[F_i^{(k)}\mid
i\in I, k\in\bbN]$ et
$
\Dist(U_\bbZ^+)=
\bbZ[E_i^{(k)}\mid
i\in I, k\in\bbN]$.
Posons $\cA=\bbZ[v,v^{-1}]$ (alg\`ebre des polyn\^omes de Laurent \`a coefficients entiers en une ind\'etermin\'ee 
$v$) et soit   maintenant $\frU_\cA$ la $\cA$-alg\`ebre enveloppante quantique    de
$G$   \`a puissances divis\'ees introduite par  Lusztig \cite{L90}.  La sp\'ecialisation 
$v\mapsto1$
fournit un isomorphisme d'anneaux
\begin{equation}
\{\frU_\cA/(K_i-1\mid i\in[1,\ell])\}
\otimes_\cA\bbZ
\simeq
\Dist(G_\bbZ)
\end{equation}
tel que  l'action $T_{i,1}''$ du groupe de tresses consid\'er\'ee par Lusztig   \cite[39.4.4]{L93}
  sur $\frU_\cA$
(et  \'egale \`a celle, $T_{\alpha_i}$,  consid\'er\'ee par  
Jantzen 
 \cite[Chap. 8]{J95},   mais diff\'erant par un signe de celle,
 $T_i$,  consid\'er\'ee par Lusztig   \cite[39.4.4]{L93}  
\cite[Thm. 3.1]{L90},   par un signe
\cite[2.1]{K00})   
induise  $s_i$ pour tout 
$i\in[1,\ell]$.

On notera que pour tous $i\in[1,\ell]$, $n\in\bbN$, on a 
\begin{align*}
\Dist(Fr)
:
E_i^{(n)}
&\mapsto
\begin{cases}
E_i^{(\frac{n}{p})}
&\text{si $p|n$}
\\
0
&\text{sinon},
\end{cases}
\\
F_i^{(n)}
&\mapsto
\begin{cases}
F_i^{(\frac{n}{p})}
&\text{si $p|n$}
\\
0
&\text{sinon},
\end{cases}
\\
\binom{H_i}{n}
&\mapsto
\begin{cases}
\binom{H_i}{\frac{n}{p}}
&\text{si $p|n$}
\\
0
&\text{sinon}.
\end{cases}
\end{align*}

\setcounter{equation}{0}
\noindent
(1.2)
Nous allons maintenant d\'efinir des scindages de Frobenius sur  $\Dist(U^\pm)$.
Ils sont induits pour $p$ ``assez grand" par ceux introduits par Lusztig  \cite[35.1.8, 35.5.2]{L93} mais requi\`erent  un travail suppl\'ementaire pour $p=2$ lorsque  $G$ est de type $\rB_n,\rC_n,\rF_4$
et pour $p=3$ lorsque $G$ est de type $\rG_2$, auxquels cas ils n\'ecessitent d'utiliser les r\'esultats de Lusztig deux fois.

\begin{prop}
Il existe des endomorphismes de  $\bbF_p$-alg\`ebres
${'\!Fr}^{\pm}$
de 
$\Dist(U^\pm)$
tels que
$E_i^{(n)}\mapsto E_i^{(np)}$
et
$F_i^{(n)}\mapsto F_i^{(np)}$,
 pour tout  $i\in[1,\ell], n\in\bbN$.

\end{prop}

\pf
Si $G$ est simplement connexe, ces endomorphismes sont obtenus \`a partir de ceux consid\'er\'es par Lusztig comme indiqu\'e pr\'ec\'edemment.
 Consid\'erons maintenant par exemple le cas $p=2$ et   $G$ est de type
$B_\ell$ avec comme racines simples
$1\Leftarrow
2
-
3
-
\cdots
-
\ell$.
Notons comme ci-dessus $U_\cA^+$
l'alg\`ebre enveloppante quantique sur
$\cA$ de type $B_\ell$.
Soit $\cA'=\cA/(v^2+1)$ et posons $U^+_{\cA'}=U^+_\cA\otimes_\cA\cA'$.
Posons   $d_1=1$ et $d_i=2$, on a alors 
 $d_s\langle\alpha_t, \alpha_s^\vee\rangle=
d_t\langle\alpha_s, \alpha_t^\vee\rangle$
pour tous $ s,t\in[1,\ell]$).  
On dispose  donc, par suite (cf.
\cite[Thm. 35.1.8]{L93}),
d'un homomorphisme de 
$\cA'$-alg\`ebres
$U^{*+}_{\cA'}\to
U^+_\cA$
tel que pour tous $ i\in[1,\ell]$, $n\in\bbN$,
\[
E_i^{*(n)}\mapsto
\begin{cases}
E_1^{(n)}
&\text{si $i=1$}
\\
E_i^{(2n)}
&\text{sinon},
\end{cases}
\]
avec
$U_{\cA'}^*$
l'alg\`ebre enveloppante quantique duale sur  $\cA'$ de type
$C_\ell$
et
$U^{*+}_{\cA'}$ sa partie positive de
g\'en\'erateurs
$E_i^{*(n)}$.
Il existe \'egalement un homomorphisme de $\cA'$-alg\`ebres  
$U^{**+}_{\cA'}\to
U^{*+}_\cA$
tel que pour tous  $i\in[1,\ell]$,
$n\in\bbN$,
\[
E_i^{**(n)}\mapsto
\begin{cases}
E_1^{*(2n)}
&\text{si $i=1$}
\\
E_i^{*(n)}
&\text{sinon}.
\end{cases}
\]
Le compos\'e de ces deux homomorphismes est donc un  homomorphisme de 
$\cA'$-alg\`ebres  
$U^{**+}_{\cA'}\to
U^{+}_\cA$
tel que pour tous   $i\in[1,\ell]$, $n\in\bbN$,
$
E_i^{**(n)}\mapsto
E_i^{(2n)}$,
ce qui donne par changement de base
$\cA'\to\bbF_2$
(avec  $v$ sp\'ecialis\'e en
$1$) un homomorphisme de 
 $\bbF_2$-alg\`ebres  
$\Dist(U^+)\to
\Dist(U^+)$
tel que pour tous  $i\in[1,\ell]$,
$n\in\bbN$,
$
E_i^{(n)}\mapsto
E_i^{(2n)}$.
De m\^eme pour  $G$ de type
$\rC_\ell$ ou
$\rF_4$
en caract\'eristique  2,
et pour  $G$ de type
$\rG_2$
en caract\'eristique  3. Finalement, sur
$\Dist(U)$ on peut transporter  
$\phi^+$
par l'anti-automorphisme
$\Omega$ de
$\Dist(G)$ (\cite[1.1]{L90}).

\setcounter{equation}{0}
\noindent
 Si maintenant, on note $B^{+}:= U^{+}T$ le sous-groupe de Borel de $G$ oppos\'e de $B$, on a : 

(1.3) {\bf Corollaire.}
{\it
Les homomorphismes de $\bbF_p$-alg\`ebres
 ${\,'\!Fr}^{\pm}$
s'\'etendent canoniquement en des homomorphismes de $\bbF_p$-alg\`ebres
 ${\,'\!Fr}^{\geq0} : \Dist(B^+)\to\Dist(B^+)$
et
$\,'\!Fr^{\leq0} : \Dist(B)\to\Dist(B)$,
tels que pour tous 
$i\in[1,\ell]$,
$n\in\bbN$,
\[
\binom{H_i}{n}\mapsto
\binom{H_i}{np}.
\]

}

\pf
Il suffit de le v\'erifier sur 
$\Dist(B^+)$ (et d'utiliser $\Omega$ une nouvelle fois pour conclure).
 On d\'efinit tout d'abord un homomorphisme 
$\bbF_p$-lin\'eaire
${\,'\!Fr}^{0}:\Dist(T)\to\Dist(T)$
par
$
\binom{H_i}{n}\mapsto
\binom{H_i}{np}$
pour tous $ i\in[1,\ell]$, $n\in\bbN$.
D'apr\`es
\cite[I.7.8.3]{J}, on a 
\begin{align*}
\binom{H_i}{n}\binom{H_i}{m}
&=
\sum_{k=0}^{\min\{n,m\}}
\frac{(n+m-k)!}{(n-k)!(m-k)!k!}\binom{H_i}{n+m-k}
\\
&=
\sum_{k=0}^{\min\{n,m\}}
\binom{n+m-k}{m}\binom{n}{k}\binom{H_i}{n+m-k}.
\end{align*}
Supposons    $n\geq m$. On a alors 
\begin{align*}
\binom{H_i}{np}\binom{H_i}{mp}
&=\sum_{k=0}^{mp}
\binom{np+mp-k}{mp}\binom{np}{k}\binom{H_i}{np+mp-k}
\\
&=
\sum_{k=0}^{m}
\binom{np+mp-kp}{mp}\binom{np}{kp}\binom{H_i}{np+mp-kp}
\\
&=
\sum_{k=0}^{m}
\binom{n+m-k}{m}\binom{n}{k}\binom{H_i}{(n+m-k)p},
\end{align*}
et donc
\[
{\,'\!Fr}^{0}(\binom{H_i}{n}\binom{H_i}{m})=
{\,'\!Fr}^{0}(\binom{H_i}{n})
{\,'\!Fr}^{0}(\binom{H_i}{m}).
\]
Comme
$\Dist(T)$
est ab\'elien et admet pour base sur   
$\bbF_p$ les
$\prod_{i\in[1,\ell]}\binom{H_i}{n_i}$,
$n_i\in\bbN$,
il en r\'esulte que 
${\,'\!Fr}^{0}$
est bien un endomorphisme de la
$\bbF_p$-alg\`ebre
$\Dist(T)$.

On remarque  ensuite que la multiplication   $\Dist(U^+)\otimes\Dist(T)\to \Dist(B^+)$ est un isomorphisme $\bbF_p$-lin\'eaire et l'on est  donc simplement ramen\'e \`a  v\'erifier, gr\^ace \`a  (1.1), que la formule de commutation entre les $E_i^{(n)}$ et les $\binom{H_j}{m}$ 
est pr\'eserv\'ee par
$\phi^+\otimes\phi^0$.
Rappelons (\cite[Lemma 26.3D]{Hum} et \cite[4.1(d)]{L89}) alors que 
\[
E_i^{(n)}
\binom{H_j}{m}
=
\binom{H_j-n\langle\alpha_i,\alpha_j^\vee\rangle}{m}
E_i^{(n)}
=
\sum_{k=0}^{m}
\binom{-n\langle\alpha_i,\alpha_j^\vee\rangle}{k}\binom{H_i}{m-k}
E_i^{(n)}.
\]
et donc
\begin{align*}
E_i^{(np)}
\binom{H_j}{mp}
&=
\sum_{k=0}^{mp}
\binom{-np\langle\alpha_i,\alpha_j^\vee\rangle}{k}\binom{H_i}{mp-k}
E_i^{(np)}
=
\sum_{k=0}^{m}
\binom{-n\langle\alpha_i,\alpha_j^\vee\rangle}{k}\binom{H_i}{(m-k)p}
E_i^{(np)},
\end{align*}
comme voulu.

\setcounter{equation}{0}
\noindent
(1.4)
On va maintenant ``recoller"
${\,'\!Fr}^{\geq0}$ et ${\,'\!Fr}^{\leq0}$
pour d\'efinir un endomorphisme de
$\bbF_p$-alg\`ebres sur tout
$\Dist(G)$  (on remarquera que  ${\,'\!Fr}^{+}\otimes{\,'\!Fr}^{0}\otimes{\,'\!Fr}^{-}$
ne convient pas). On introduit pour se faire  une mesure invariante 
 sur $\Dist(T_1)$ (avec $T_1$ le noyau de Frobenius de $T$) en posant
\[
\mu_0=
\prod_i(\sum_j^{p-1}(-1)^j\binom{H_i}{j})=
\prod_i\binom{H_i-1}{p-1}
=
\prod_i(1-H_i^{p-1})
\]
C'est une involution de  
$\otimes_{i\in[1,\ell]}\{(\Dist(T_{\alpha_i,1})^{T_{\alpha_i,1}}\}=
\Dist(T_1)^{T_1}$ 
qui commute  
avec tous les 
$E_i^{(n)}$
et 
$F_i^{(n)}$,
$i\in[1,\ell]$,
$n\in\bbN$ et l'on pose  
\begin{align*}
\phi^{\leq0}&={\,'\!Fr}^{\leq0}(?)\mu_0
:
\Dist(B)\to\Dist(G),
\\
\phi^{\geq0}&={\,'\!Fr}^{\geq0}(?)\mu_0
:
\Dist(B^+)\to\Dist(G),
\\
\phi^{-}&={\,'\!Fr}^{\leq0}|_{\Dist(U)}(?)\mu_0
:
\Dist(U)\to\Dist(G),
\\
\phi^{+}&={\,'\!Fr}^{\geq0}|_{\Dist(U^+)}(?)\mu_0
:
\Dist(U^+)\to\Dist(G),
\\
\phi^{0}&={\,'\!Fr}^{\leq0}|_{\Dist(T)}(?)\mu_0
={\,'\!Fr}^{\geq0}|_{\Dist(T)}(?)\mu_0 :
\Dist(T)\to\Dist(G).
\end{align*}
Comme ${\,'\!Fr}^{\leq0}$
et ${\,'\!Fr}^{\geq0}$ sont toutes deux des applications multiplicatives,
les applications $\phi^{\lesseqgtr0}$,
$\phi^{\pm}$
et $\phi^0$
le sont aussi.
On va voir que l'on peut d\'efinir  
une application  $\bbF_p$-lin\'eaire
$\phi : \Dist(G)\to\Dist(G)$
par
$FHE\mapsto\phi^-(F)\phi^0(H)\phi^+(E)$
pour tous $F\in\Dist(U),
H\in\Dist(T),
E\in\Dist(U^+)$.

\begin{thm}
L'application   $\phi$ est multiplicative et induit  donc un homomorphisme de $\bbF_p$-alg\`ebres $\phi : \Dist(G)\to \im \,\phi$
(avec
$\mu_0$ comme unit\'e de $\im\phi$). On a, de plus,     $\Dist(Fr)\circ\phi=\id_{\Dist(G)}$.
 \end{thm}

\pf
Montrons que pour tous 
$F'\in\Dist(U),
H'\in\Dist(T),
E'\in\Dist(U^+)$, on a 
\[
\phi(FHEF'H'E')=
\phi(FHE)\phi(F'H'E').
\]
Posons
$E_i=X_{\alpha_i}$
et
$F_i=X_{-\alpha_i}$
pour tout $\alpha_i\in R^s$.
On peut supposer par (1.1) que
$E$ (resp. $F'$)
est un produit de $E_i^{(n_i)}$,
$n_i\in\bbN$
(resp. de
$F_i^{(m_i)}$,
$m_i\in\bbN$).
On \'ecrit alors
$E=\prod_s
E_{i_s}^{(n_{i_s})}$,
$i_s\in[1,\ell]$,
$n_{i_s}\in\bbN$,
et
$F=\prod_t
F_{i_t}^{(m_{i_t})}$,
$i_t\in[1,\ell]$,
$m_{i_t}\in\bbN$ et l'on va raisonner par r\'ecurrence sur
$\sum_s
n_{i_s}+\sum_t
m_{i_t}$.

Si $\sum_sn_{i_s}=0$,
\begin{align*}
\phi(FHF'H'E')
&=
\phi^{\leq0}(FHF'H')\phi^+(E')
\quad\text{par d\'efinition car
$FHF'H'\in\Dist(B)$}
\\
&=
\phi^-(F)\phi^0(H)
\phi^-(F')\phi^0(H')\phi^+(E')
\quad\text{car
$\phi^{\leq0}$
est multiplicative}
\\
&=
\phi(FH)
\phi(F'H'E').
\end{align*}
De m\^eme, si
$\sum_tm_{i_t}=0$,
\begin{align*}
\phi(FHEH'E')
&=
\phi^-(F)\phi^{\geq0}(HEH'E')
\quad\text{par d\'efinition car
$HEH'E'\in\Dist(B^+)$}
\\
&=
\phi^-(F)\phi^0(H)
\phi^+(E)\phi^0(H')\phi^+(E')
\quad\text{car
$\phi^{\geq0}$
est multiplicative}
\\
&=
\phi(FHE)
\phi(H'E').
\end{align*}

Supposons maintenant que
$E=E_i^{(n_i)}$ avec
$n_i>0$ et
$F'=F'_1F_i^{(m_i)}F'_2$
avec
$F'_1$ ne contenant pas de facteurs contenant des puissances divis\'ees des $F_i$.
Si $m_i>0$,
\begin{align*}
\phi(
&FHE_i^{(n_i)}F'H'E')
=
\phi(FHE_i^{(n_i)}F'_1F_i^{(m_i)}F'_2H'E')
\\
&=
\phi(FH
F_1'E_i^{(n_i)}F_i^{(m_i)}F'_2H'E')
\\
&=
\phi(FHF'_1\sum_{k=0}^{\min\{n_i,m_i\}}
F_i^{(m_i-k)}
\binom{H_i-n_i-m_i+2k}{k}
E_i^{(n_i-k)}F'_2H'E')
\\
&=
\sum_{k=0}^{\min\{n_i,m_i\}}
\phi(FHF'_1F_i^{(m_i-k)}
\binom{H_i-n_i-m_i+2k}{k}
E_i^{(n_i-k)})\phi(
F'_2H'E')
\\
&\hspace{1cm}\text{par hypoth\`ese de r\'ecurrence car le nombre de facteurs dans $F'_2$
est moindre que celui dans  
$F'$}
\\
&=
\sum_{k=0}^{\min\{n_i,m_i\}}
\phi^{\leq0}(FHF'_1)
\phi^-(F_i^{(m_i-k)})
\phi^0(\binom{H_i-n_i-m_i+2k}{k})
\phi^+(E_i^{(n_i-k)})
\phi(F'_2H'E')
\\
&=
\phi^{\leq0}(FHF'_1)
\phi(E_i^{(n_i)}F_i^{(m_i)})
\phi(F'_2H'E')
\quad\text{par d\'efinition}
\\
&=
\phi^{\leq0}(FHF'_1)
\phi^+(E_i^{(n_i)})\phi^-(F_i^{(m_i)})
\phi(F'_2H'E')
\quad\text{by \cite{G}}
\\
&=
\phi^{\leq0}(FH)\phi^-(F'_1)
\phi^+(E_i^{(n_i)})\phi^-(F_i^{(m_i)})
\phi^-(F'_2)\phi^{\geq0}(H'E')
\\
&=
\phi^{\leq0}(FH)\phi^+(E_i^{(n_i)})
\phi^-(F'_1)
\phi^-(F_i^{(m_i)})
\phi^-(F'_2)\phi^{\geq0}(H'E')
\\
&=
\phi^{\leq0}(FH)\phi^+(E_i^{(n_i)})
\phi^-(F')\phi^{\geq0}(H'E')
\\
&=
\phi(FHE_i^{(n_i)})
\phi(F'H'E')
\quad\text{par d\'efinition}.
\end{align*}
Si aucun facteur avec des puissances divis\'ees de $F_i$ n'apparait dans   
$F'$,
\begin{align*}
\phi(
&FHE_i^{(n_i)}F'H'E')
=
\phi(FHF'E_i^{(n_i)}H'E')
\\
&=
\phi^{\leq0}(FHF')\phi^{\geq0}(E_i^{(n_i)}H'E')
\quad\text{par d\'efinition}
\\
&=
\phi^{\leq0}(FH)\phi^-(F')
\phi^+(E_i^{(n_i)})\phi^{\geq0}(H'E')
=
\phi^{\leq0}(FH)\phi^+(E_i^{(n_i)})\phi^-(F')
\phi^{\geq0}(H'E')
\\
&=
\phi(FHE_i^{(n_i)})\phi(F'H'E')
\quad\text{par d\'efinition de nouveau}.
\end{align*}
Finalement,
si $E=E''E_i^{(n_i)}$ avec $n_i>0$,
\begin{align*}
\phi(
&FHE''E_i^{(n_i)}F'H'E')
\\
&=
\phi(FHE'')\phi(E_i^{(n_i)}F'H'E')
\\
&\hspace{1cm}
\text{par r\'ecurrence car le nombre de facteurs dans
 $E''$ est moindre que celui dans
$E$}
\\
&=
\phi(FHE'')\phi(E_i^{(n_i)})\phi(F'H'E')
\quad\text{gr\^ace au cas 
$E=E_i^{(n_i)}$
plus haut}
\\
&=
\phi^{\leq0}(FH)\phi^+(E'')
\phi^+(E_i^{(n_i)})\phi(F'H'E')
\quad\text{par d\'efinition}
\\
&=
\phi^{\leq0}(FH)\phi^+(E''
E_i^{(n_i)})\phi(F'H'E')
\\
&=
\phi(FHE''
E_i^{(n_i)})\phi(F'H'E')
\quad\text{par d\'efinition de nouveau},
\end{align*}
ce qui conclut.

\setcounter{equation}{0}
\noindent
(1.5)
{\it Remarque.} Nous remercions T. Tanisaki de nous avoir signal\'e que lorsque  l'on s'int\'eresse seulement \`a l'alg\`ebre quantique ``modifi\'ee", alors l'existence d'un scindage de Frobenius (qui est simplement le recollement de celui d\'efini sur les parties nilpotentes et torique)  est d\'ej\`a connu  (\cite[Prop. 3.4] {Mc}).

\begin{center}
{\bf
$2^\circ$
Contractions}
\end{center}

\setcounter{equation}{0}
\noindent
(2.1)
Si  $M$ est un $G$-module,
$M^{[1]}$ d\'esignera le $G$-module 
twist\'e par Frobenius, ou de mani\`ere \'equivalente
le $\bbF_p$-espace vectoriel
$M$ muni de l'action de 
$\Dist(G)$ via
$\Dist(Fr)$.
Soit $\chi_\lambda$
l'homomorphisme de $\bbF_p$-alg\`ebres  
de
$\Dist(T)$ vers $\bbF_p$
d\'efini par $\chi_\lambda(\displaystyle\binom{H_i}{n})=
\binom{\langle\lambda,\alpha_i^\vee\rangle}{n}
\ $   pour tous $i\in[1,\ell]$, $n\in\bbN$.
On a
$\chi_\lambda\circ\Dist(Fr)|_{\Dist(T)}
=
\chi_{p\lambda}$.

La mesure $\mu_0$ \'etant  une involution,
tout $M\in\cC_G$ admet une d\'ecomposition
$M=\mu_0M\oplus(1-\mu_0)M$.
Comme on a muni ${\text{im}}\, \phi$ de
$\mu_0$  comme unit\'e, on peut d\'efinir une structure de
 $G$-module sur $\mu_0M$ 
via
$\phi$, annulant
$(1-\mu_0)M$.
Nous noterons cette nouvelle  structure sur
$\mu_0M$ par $M^\phi$, et \'ecrirons
$x\bullet m=\phi(x)m$
avec
$x\in\Dist(G)$ et $m\in\mu_0M$.
Nous noterons l'action correspondante de   $G$ induite par
$\phi$ par la m\^eme lettre.
De m\^eme pour $M\in\cC_B, \cC_T$.
On appellera  $M^\phi$ la contraction par  Frobenius de $M$.

Soit
$\Lambda_1=\{\lambda\in\Lambda\mid  \forall i\in[1,\ell]\ ,
\langle \lambda,\alpha_i^\vee\rangle\in[0,p[\}$.
Pour $\lambda\in\Lambda$
nous \'ecrirons
$\lambda=\lambda^0+p\lambda^1$
avec
$\lambda^0\in\Lambda_1$
et
$\lambda^1\in\Lambda$.

\begin{lem}
Pour tout $\lambda\in\Lambda$,
\[
\chi_\lambda\circ\phi|_{\Dist(T)}
=
\begin{cases}
\chi_{\lambda^1}
&\text{si $\lambda\in p\Lambda$},
\\
0
&\text{sinon}.
\end{cases}
\]
En particulier, pour tout
$M\in\cC_T$,
\[
M^\phi=
\coprod_{\lambda\in\Lambda}M_{p\lambda}.
\]

\end{lem}

\pf
On a 
\begin{align*}
\chi_\lambda(\mu_0)
&=
\chi_\lambda(\prod_{i=1}^\ell
(1-H_i^{p-1}))
=
\prod_{i=1}^\ell(1-\chi_\lambda(H_i)^{p-1})
=
\prod_{i=1}^\ell(1-\langle\lambda,
\alpha_i\rangle^{p-1})
\\
&=
\begin{cases}
1
&\text{si $p|\langle\lambda,\alpha_i^\vee\rangle$
pour tout $i$}
\\
0
&\text{sinon}.
\end{cases}
\end{align*}
Il en d\'ecoule alors que pour tout
$m\in\bbN$,
\begin{align*}
\chi_\lambda\circ\phi
\binom{H_i}{m}
&=
\chi_\lambda
\binom{H_i}{mp}
\chi_\lambda(\mu_0)
=
\begin{cases}
\binom{\langle\lambda,\alpha_i^\vee\rangle}{mp}=
\binom{p\langle\lambda^1,\alpha_i^\vee\rangle}{mp}=
\binom{\langle\lambda^1,\alpha_i^\vee\rangle}{m}
&\text{si $\lambda\in p\Lambda$
}
\\
0
&\text{sinon}
\end{cases}
\end{align*}
et comme l'on a
$\chi_{\lambda^1}
\binom{H_i}{m}
=
\binom{\langle\lambda^1,\alpha_i^\vee\rangle}{m}$, le lemme en r\'esulte.

\noindent
(2.2)
Soient $W$ le groupe de Weyl de $G$ et $M\in\cC_G$. Si $w\in W$, un repr\'esentant de $w$ dans $\rN_G(T)$ permute  les espaces de poids  de
$M$.

\begin{lem}
Si $m\in M_{p\lambda}$,
$\lambda\in\Lambda$, alors
$w\bullet m=\phi(w)m\in
M_{pw\lambda}$.

\end{lem}

\pf
Pour tout $h\in\Dist(T)$,
$h\bullet(w\bullet m)=w\bullet\Ad(w^{-1})(h)\bullet m=
\chi_{w\lambda}(h)w\bullet m$. L'assertion r\'esulte alors de   (2.1).

\noindent

\noindent
(2.3)
{\it Exemples.}
Soit
$M\in\cC_G$.

(i) Le $G$-module
$(M^{[1]})^\phi$ s'identifie \`a $M$.

(ii)
Soient $M\in\cC_G$ et $m\in M$.
Si $\Delta_M : M\to M\otimes_{\bbF_p}\bbF_p[t]$
est le comorphisme pour la $U_{\alpha_i}$-structure sur $M$, on a 
$\Delta_M(m)=\sum_{r\in\bbN}E_i^{(r)}m\otimes t^r$ par \cite[I.7.2]{J}.
Alors
$\Delta_{M^\phi}(m)=\sum_{s\in\bbN}E_i^{(ps)}m\otimes t^s$.
On en d\'eduit, pour $A$ une $\bbF_p$-alg\`ebre commutative, 
que
$\phi(x_{\alpha_i}(a))m=
\sum_{s\in\bbN}a^sE_i^{(ps)}m$
pour tout $a\in A$.

(iii)
Soit $G=\SL_2$ et identifions $\Lambda$ avec
$\bbZ$. Notons, pour $n\in\bbN$,  $\nabla(n)$ le  $G$-module standard induit \`a partir  du $B$-module de de dimension 1 de poids $n$  
et  $L(n)$ son socle simple. 

Sur $\bbF_2$
on a
$\nabla(2)^\phi\simeq
L(1)\oplus L(0)\simeq
\Delta(2)^\phi$,
$(\nabla(2)^\phi)^{[1]}\simeq
L(2)\oplus L(0)\simeq
(\Delta(2)^\phi)^{[1]}$,
$(\St\otimes M^{[1]})^\phi=0$.

Sur $\bbF_3$ 
on a
$\nabla(3)^\phi
=L(1)=\Delta(3)^\phi$,
$(\nabla(3)^\phi)^{[1]}=L(3)=(\Delta(3)^\phi)^{[1]}$,
et donc
$\nabla(3)^\phi\simeq
\nabla(3)/L(3)$,
$(\nabla(3)^\phi)^{[1]}<\nabla(3)$
alors que
$\Delta(3)^\phi<\Delta(3)$,
$(\Delta(3)^\phi)^{[1]}\simeq
\Delta(3)/L(1)$.

Maintenant, sur $\bbF_p$ avec
$p\geq5$, 
$\nabla(p)^\phi
=L(1)=\Delta(p)^\phi$
n'est ni un facteur de Jordan-H\"older de
$\nabla(p)$
ni de  $\Delta(p)$
alors que
$(\nabla(p)^\phi)^{[1]}=L(p)=(\Delta(p)^\phi)^{[1]}$ est \`a la fois le socle de
$\nabla(p)$
et le module de t\^ete de 
$\Delta(p)$.

Sur $\bbF_p$ avec
$p$ impair, $(\St\otimes M^{[1]})^\phi\simeq
M$.

\begin{center}
{\bf
$3^\circ$
La th\'eorie de Kumar-Littelmann }
\end{center}

Son principe est le suivant : si
$X$ est une vari\'et\'e  projective,
on peut, par d\'efinition \'ecrire
$X=\Proj(A)$ avec $A$ un anneau gradu\'e,
et \'etudier la g\'eom\'etrie de $X$ \`a partir des propri\'et\'es de  $A$.
En particulier, soient
$\lambda\in \Lambda$ un poids dominant,
$I=\{\alpha\in R^s\mid\langle\lambda,\alpha^\vee\rangle=0\}$,
et $P=P_{I}$ le sous-groupe parabolique de $G$ standard associ\'e \`a $I$.
Si
$\nabla
(\lambda)$ d\'esigne le $G$-module induit \`a partir du 
  $P$-module de dimension un de poids 
$\lambda$ et si
$v_-$ est un vecteur de plus bas poids de  
$\nabla
(\lambda)^*$,
il existe une immersion ferm\'ee  
$G/P\to\bbP(\nabla
(\lambda)^*)$
via
$g\mapsto[gv_-]$.
On a alors
$G/P\simeq\Proj(
A)$
avec
$A
=
\coprod_{n\in\bbN}\nabla
(n\lambda)$
vu comme anneau gradu\'e via le  ``cup-produit''.
Notre contraction par Frobenius est adapt\'ee \`a la construction d'un scindage    $G_\Bbbk$-``semi-invariant"  du comorphisme de Frobenius  $\cO_{G_\Bbbk/P_\Bbbk}\to F_*\cO_{G_\Bbbk/P_\Bbbk}$,
comme d\'ecrit  dans  \cite{KuLi} avec la contraction de  Littelmann 
\cite{Li}. Comme les d\'etails sont essentiellement les m\^emes que dans   \cite{KuLi}, seulement rendus plus simples ici et l\`a par l'absence de techniques quantiques, nous ne ferons qu'esquisser  la construction.  

\noindent

\noindent
(3.1)
Rappelons n\'eanmmoins tout d'abord le formalisme des foncteurs d'induction tel  que formul\'e  par  Zuckermann dans le cas classique   \cite{EW} 
 et par Andersen, Polo et Wen dans le cas quantique
\cite{APW91}.
Soit $G\Mod$ la cat\'egorie des $G$-modules rationnels.
Elle est \'equivalente \`a la cat\'egorie 
$\cC_G$
des
$\Dist(G)$-modules int\'egrables,
i.e.  des
$\Dist(G)$-modules localement finis
$M$
qui admettent une d\'ecomposition en espace de poids relativement \`a l'action de  $T$  \cite[II.1.20]{J}/\cite{CPS}.
Rappelons aussi  
 \cite[I.7.16]{J} 
que pour tous
$  M, M'\in\cC_G$, on a 
$G\Mod(M,M')=\Dist(G)\Mod(M,M')$.
De m\^eme avec tout sous-groupe parabolique standard  
 $P$ de $G$ contenant  $B$ et aussi pour  $T$:
$P\Mod\simeq
\cC_P, T\Mod\simeq\cC_T$,
ainsi que  pour les morphismes.

Soit $\ind_P^G=
\Sch_{\bbF_p}(G, ?)^P$ le foncteur d'induction de
 $P\Mod$ \`a  $G\Mod$ avec comme morphisme d'\'evaluation   $\ev$, l'\'evaluation en    l'\'el\'ement neutre
$e$ de $G$.
Via l'\'equivalence de cat\'egories rappel\'ee ci-dessus, ce foncteur est \'equivalent \`a l'adjoint \`a droite  $\cF_P^G$
du foncteur de restriction de  $\cC_G$ vers $\cC_P$.
Plus pr\'ecis\'ement, si $M\in\cC_P$, on a 
$\cF_P^GM=
\{f\in
\coprod_{\lambda\in\Lambda}
\Dist(P)\Mod(\Dist(G), M)_\lambda\mid
E_i^{(n)}f=0=F_i^{(n)}f\  
{\text{pour tous}}\, i\in[1,\ell] \,,\,  n  >>0\}$,
avec
la $\Dist(G)$-action sur
$\Dist(P)\Mod(\Dist(G), M)$
donn\'ee par
$(xf)(y)=f(yx)$,
$x, y\in\Dist(G)$,
et pour tout
$\lambda\in\Lambda$,
$\Dist(P)\Mod(\Dist(G), M)_\lambda:=
\{f\in\Dist(P)\Mod(\Dist(G), M)\mid
\displaystyle\binom{H_i}{n}f=
\binom{\langle\lambda,\alpha_i^\vee\rangle}{n}f\ 
{\text{pour tous}}\, i\in[1,\ell]\ ,  n\in\bbN\}$. L'ensemble 
$\cF_P^GM$ est muni d'une action  
$P$-lin\'eaire, i.e., d'une application $\Dist(P)$-lin\'eaire
$\ev_M : \cF_P^GM\to M$
via
$f\mapsto
f(1)$.
Donc, pour tout
$Q\in\cC_G$,
$
\cC_P(Q,M)\simeq\cC_G(Q,\cF_P^GM)
$ via
$\psi\mapsto
\lq\lq q\mapsto\psi(?q)"$
pour tout $q\in Q$
d'inverse
$\eta\mapsto\ev_M\circ\eta$.
En particulier, on a un isomorphisme
$\ind_P^GM\simeq\cF_P^GM$
via
$h\mapsto(?h)(e)$
d'inverse
$(?^{-1}f)(1)\gets\!\shortmid f$.

\setcounter{equation}{0}
\noindent
(3.2)
Soient 
$\bfR^\bullet\cF_P^G$ les foncteurs d\'eriv\'es droits de  
$\cF_P^G$. On a donc
$\bfR^\bullet\cF_P^G(M)\simeq
\rH^\bullet(G/P, \cL(M))$
avec
$\cL(M)$ le $\cO_{G/P}$-module
associ\'e au $P$-module $M$.
Soit
$\Lambda_P$ le groupe des caract\`eres de $P$. Si
$R^s_P$ d\'esigne l'ensemble des racines simples du sous-groupe de Levi  standard de $P$, on a
$\Lambda_P=\{\lambda\in\Lambda\mid
\langle\lambda,\alpha^\vee\rangle=0\}$ et l'on pose 
$\Lambda_P^+=\Lambda_P\cap\Lambda^+$.
Posons aussi $I_P=\{i\in[1,\ell]\mid\alpha_i\in R_P^s\}$, et soit $U_P^+$ le radical unipotent du sous-groupe parabolique oppos\'e,
i.e., le sous-groupe engendr\'e par les sous-groupes radiciels de    $G$ associ\'es aux racines
$\alpha\in R^+\setminus\sum_{i\in I_P}\bbN\alpha_i$.
Si $\lambda\in\Lambda_P^+$,
le $G$-module $\ind_P^G(\lambda)\simeq\cF_P^G\lambda$ est isomorphe \`a  
$\ind_B^G(\lambda)\simeq\cF_B^G\lambda$, que nous noterons  
$\nabla(\lambda)$.

\noindent
{\bf
Proposition}
(cf. \cite[Lem. 2.2, Th. 2.3]{KuLi}).
{\it
Pour tout $M\in\cC_B$, il existe un morphisme $G$-lin\'eaire (ie. $\Dist(G)$-lin\'eaire)
$\Phi_M : (\cF_P^GM)^{[1]}\to
\cF_P^G(M^{[1]})$
tel que
$(\Phi_Mf)(x)=f(\Dist(Fr)(x))$
pour tous $  f\in\cF_B^GM$,
$ x\in\Dist(G)$. Ce morphisme est fonctoriel en
$M$ et induit des applications 
$G$-lin\'eaires  
$\bfR^\bullet\Phi_M : (\bfR^\bullet\cF_B^GM)^{[1]}\to
\bfR^\bullet\cF_B^G(M^{[1]})$.
En particulier,
lorsque $M=\lambda$,
l'application $\lambda\in\Lambda$,
$\Phi_\lambda : (\cF_B^G\lambda)^{[1]}\to
\cF_B^G(p\lambda)$ n'est autre que l'application d'\'el\'evation \`a la puissance $p$, via le cup-produit.

}

\setcounter{equation}{0}
\noindent
(3.3)
\noindent
{\bf
Proposition}
(cf. \cite[Lem. 3.2]{KuLi}). 
{\it
Pour  tout $ M\in\cC_T$, il existe un  morphisme  $B$-lin\'eaire (i.e. $\Dist(B)$-lin\'eaire), fonctoriel en $M$,
$\Phi_{B,M} : (\cF_T^BM)^{\phi}\to
\cF_T^B(M^{\phi})$
tel  que
$(\Phi_{B,M}f)(x)=f(\phi(x))$ pour tous 
$f\in\cF_T^BM$,
$x\in\Dist(B)$.

}

\setcounter{equation}{0}
\noindent
(3.4)
\noindent
{\bf
Proposition}
(cf. \cite[Prop. 3.3, Th. 3.4]{KuLi}). 
{\it
Pour tout $M\in\cC_B$, il existe un  morphisme  $G$-lin\'eaire (i.e. $\Dist(G)$-lin\'eaire) fonctoriel en $M$
$\Psi_M : (\cF_B^GM)^{\phi}\to
\cF_B^G(M^{\phi})$
tel que
$(\Psi_Mf)(x)=f(\phi(x))$
pour tous $f\in\cF_B^GM$,
$x\in\Dist(G)$,
  induisant des applications   
$G$-lin\'eaires  naturelles
$\bfR^\bullet\Psi_M : (\bfR^\bullet\cF_B^GM)^{\phi}\to
\bfR^\bullet\cF_B^G(M^{\phi})$.

}

\setcounter{equation}{0}
\noindent
(3.5)
{\it Remarques.}
(i)
Si $\lambda\in\Lambda\setminus p\Lambda$,
alors $\lambda^\phi=0$
par d\'efinition, et donc
$\bfR^\bullet\Psi_\lambda :
(\bfR^\bullet\cF_B^G\lambda)^\phi\to
\bfR^\bullet\cF_B^G(\lambda^\phi)$
est l'application nulle.

(ii)
(cf. \cite[Lem. 4.13]{KuLi})
Pour tous $\lambda\in\Lambda\setminus p\Lambda$,
$f\in\cF_B^G(M^{[1]})_\lambda$, on a
$\Psi_Mf=0$.

\pf
Comme
$\Psi_Mf\in\Dist(B)\Mod(\Dist(G),M)$,
il suffit de v\'erifier que
$\Psi_Mf$ annule $\Dist(U^+)$.
Soit
$y\in\Dist(U^+)$.
On peut supposer que
$y=E_{i_1}^{(m_{i_1})}\dots
E_{i_s}^{(m_{i_s})}$ avec $i_k\in[1,\ell]$, $m_{i_k}\in\bbN$.
Par hypoth\`ese, il existe
$j\in[1,\ell]$ tel que
$\langle\lambda,\alpha_j^\vee\rangle=
\lambda_j^0+p\lambda_j^1$
avec
$\lambda_j^0\in]0,p[$.
On a alors
$f(\phi(y)\binom{H_j}{\lambda_j^0}))
=
(\binom{H_j}{\lambda_j^0})f)(\phi(y))=
\binom{\lambda_j^0+p\lambda_j^1}{\lambda_j^0}f(\phi(y))=
f(\phi(y))=
(\Psi_Mf)(y)$.
D'autre part,
\begin{align*}
f(\phi(y)\binom{H_j}{\lambda_j^0})
&=
f(E_{i_1}^{(m_{i_1})}\dots
E_{i_s}^{(m_{i_s})}\mu_0\binom{H_j}{\lambda_j^0})
\\
&=
f(\binom{H_j-\sum_kpm_{i_k}\langle\alpha_{i_k},\alpha_j^\vee\rangle}{\lambda_j^0}E_{i_1}^{(m_{i_1})}\dots
E_{i_s}^{(m_{i_s})}\mu_0)
\quad\text{par \cite[Lem.26.3D]{Hum}}
\\
&=
\binom{H_j-\sum_kpm_{i_k}\langle\alpha_{i_k},\alpha_j^\vee\rangle}{\lambda_j^0}\bullet
f(\phi(y))
\\
&=
\sum_{r=0}^{\lambda_j^0}\binom{-\sum_kpm_{i_k}\langle\alpha_{i_k},\alpha_j^\vee\rangle}{r}\binom{H_j}{\lambda_j^0-r}\bullet
f(\phi(y))
\quad\text{par \cite[Prop. 2.1.1]{G}}
\\
&=
\sum_{r=0}^{\lambda_j^0-1}\binom{-\sum_kpm_{i_k}\langle\alpha_{i_k},\alpha_j^\vee\rangle}{r}\binom{H_j}{\lambda_j^0-r}\bullet
f(\phi(y))
\\
&\hspace{4cm}
\text{car
$\binom{-\sum_kpm_{i_k}\langle\alpha_{i_k},\alpha_j^\vee\rangle}{\lambda_j^0}=0$
par  \cite[preuve de Prop. 2.1.1]{G}}
\\
&=
\sum_{r=0}^{\lambda_j^0-1}\binom{-\sum_kpm_{i_k}\langle\alpha_{i_k},\alpha_j^\vee\rangle}{r}
(\Dist(Fr)\binom{H_j}{\lambda_j^0-r})
f(\phi(y))
=0.
\end{align*}
Il en d\'ecoule que
$(\Psi_Mf)(y)=0$.

\setcounter{equation}{0}
\noindent
(3.6)
{\bf Th\'eor\`eme}
(cf. \cite[Cor. 3.9]{KuLi}). 
{\it
Pour tout $M\in\cC_B$, il existe un diagramme commutatif de $G$-modules
\[
\xymatrix@R2ex{
((\bfR^\bullet\cF_B^GM)^{[1]})^\phi
\ar@{-}[d]_\sim
\ar[rr]^-{(\bfR^\bullet\Phi_M)^\phi}
&&
(\bfR^\bullet\cF_B^G(M^{[1]}))^\phi
\ar[dd]^-{\bfR^\bullet\Psi_{M^{[1]}}}
\\
\bfR^\bullet\cF_B^GM
\ar[rrdd]_{\id_{\bfR^\bullet\cF_B^GM}}&&
\\
&&
\bfR^\bullet\cF_B^G((M^{[1]})^\phi)
\ar@{-}[d]^-\sim
\\
&&
\bfR^\bullet\cF_B^GM.
}
\]

}

\setcounter{equation}{0}
\noindent
(3.7)
Posons $\rho=\displaystyle\frac{1}{2}\sum_{\alpha\in R^+}\alpha$ et $E^+=
\prod_{\alpha\in R^+}E_\alpha^{(p-1)} \in \Dist(G)$.
Comme $\dim\Dist(U^+_1)_{2(p-1)\rho}=1$,
$E^+$ est, \`a une constante pr\`es, dans 
$\Bbbk^\times$ ind\'ependant du choix de l'ordre dans les racines positives figurant dans le produit.  

\begin{lem}
Soit $i\in[1,\ell]$.
Pour tout $r\in\bbN$, il existe
$z_{\nu s}\in\Dist(T_1)$ (avec
$\nu\in[0,p[^{R^+}$, 
$s\in[0,r[$),
tel  que  $\chi_{2(p-1)\rho}(z_\nu)=0$ et
\[
E^+F_i^{(rp)}\in
F_i^{(rp)}
E^++
\underset{\substack{
s\in]0,rp[
\\
p\not\,| \ s}}{\sum}
F_i^{(s)}\Dist(G)+
\sum_{\nu\in[0,p[^{R^+}}\sum_{s=0}^{r-1}F_i^{(sp)}z_{\nu s}\Dist(G).
\]

\end{lem}

\pf
Nous proc\'ederons par r\'ecurrence sur  $r$. Montrons donc tout d'abord le cas $r=1$. Comme $G_1\triangleleft G$, 
l'action adjointe de $F_i^{(p)}$ sur $\Dist(G)$ laisse
$\Dist(G_1)$-invariant
\cite[3.4.13]{Ta}. 
Comme
$\Delta(F_i)=1\otimes F_i+F_i\otimes1$, on a 
$\Delta(F_i^{(p)})=
\sum_{k=0}^p
F_i^{(k)}\otimes F_i^{(p-k)}$,
et donc
\begin{align*}
\Dist(G_1)
&\ni
\ad(F_i^{(p)})(E^+)=
\sum_{k=0}^p
F_i^{(k)}E^+(-1)^{p-k}F_i^{(p-k)}
\\
&=
-E^+F_i^{(p)}+F_i^{(p)}E^+
+
\sum_{k=1}^{p-1}
F_i^{(k)}E^+(-1)^{p-k}F_i^{(p-k)}.
\end{align*}
Comme $-\alpha_i$ est l'unique racine n\'egative impliqu\'ee dans l'expression de  $\ad(F_i^{(p)})(E^+)$ dans
$\Dist(G_1)=\Dist(U_1)\Dist(T_1)\Dist(U^+_1)$
\cite[Lem.26.3C]{Hum}, on peut \'ecrire
\begin{equation}
E^+F_i^{(p)}=F_i^{(p)}E^+
+
\sum_{j=1}^{p-1}F_i^{(j)}x_j+\sum_{\nu\in[0,p[^{R^+}}
z_\nu E^{(\nu)}
\quad\text{with}\quad
E^{(\nu)}=\prod_{\alpha\in R^+}E^{(\nu_\alpha)}
\end{equation}
avec  $x_j\in\Dist(B_1^+)$
$z_\nu\in\Dist(T_1)$.
Montrons que  $\chi_{2(p-1)\rho}(z_\nu)=0$ pour tout 
$ \nu$ suivant en cela
\cite[Lem. 4.5]{KuLi}.
Rappelons  
\cite[1.1 (d2)]{L90} 
qu'il existe un anti-automorphisme $\tau$
de
$\Dist(G)$
fixant les $E_j$ et les $F_j$ 
et tel que
$\binom{H_j}{m}\mapsto
\binom{-H_j}{m}$
pour tous $ m\in\bbN$, 
$ j\in[1,\ell]$.
Appliquant $\tau$ \`a (1), on obtient
\[
F_i^{(p)}\tau(E^+)=\tau(E^+)F_i^{(p)}
+
\sum_{j=1}^{p-1}\tau(x_j)\tau(F_i^{(j)})+\sum_{\nu\in[0,p[^{R^+}}
\tau(E^{(\nu)})\tau(z_\nu).
\]
Si
$v_-\in\Delta(2(p-1)\rho)_{-2(p-1)\rho}\setminus0$,
alors
$F_i^{(p)}\tau(E^+)v_-=
\tau(E^{(\nu)})\tau(z_\nu)v_-$.
Montrons maintenant que
\begin{equation}
F_i^{(p)}\tau(E^+)v_-=0.
\end{equation}
Il suffit pour cela de montrer que
$F_i^{(p)}E^+v_-=0$.
Supposons que cela ne soit pas le cas.
Alors, comme le poids de
$F_i^{(p)}E^+v_-$ est $-p\alpha_i$,
par la th\'eorie  $\SL_2$, on doit avoir  
\begin{align*}
0&\ne
E_i^{(p)}E^+v_-
\\
&=
E^+E_i^{(p)}v_-
\quad\text{par
\cite[6.1(v)]{X}/\cite[Lem.26.C]{Hum}}
\\
&=
(\prod_{\alpha\in R^+\setminus\alpha_i}
E_\alpha^{(p-1)})E_i^{(p-1)}E_i^{(p)}v_-
=
(\prod_{\alpha\in R^+\setminus\alpha_i}
E_\alpha^{(p-1)})E_i^{(2p-1)}v_-,
\end{align*}
et donc
$E_i^{(2p-1)}v_-\in\Delta(2(p-1)\rho)_{-2(p-1)\rho+(2p-1)\alpha_i}\setminus0$.
Mais alors
$s_i(-2(p-1)\rho+(2p-1)\alpha_i)=
-2(p-1)\rho+2(p-1)\alpha_i-(2p-1)\alpha_i=
-2(p-1)\rho-\alpha_i$
devrait aussi \^etre un poids de
$\Delta(2(p-1)\rho)$,
ce qui est absurde.

Il s'ensuit que
$\tau(E^{(\nu)})\tau(z_\nu)v_-=0$.
Via le cup-produit
$\St\otimes\St\to\nabla(2(p-1)\rho)$, on a 
$v_+\otimes v_+\mapsto
v_{++}$
avec
$v_+$
(resp. $v_{++}$) un vecteur de plus haut poids de $\St$
(resp. $\nabla(2(p-1)\rho)$).
En dualisant, on obtient une application  
 $G$-lin\'eaire  
$\Delta(2(p-1)\rho)\to\St\otimes\St$
envoyant
$v_-$
vers
$v_-'\otimes v_-'$
avec
$v_-'$ un vecteur de plus bas poids de
$\St$.
Alors
$E^+v_-$
est envoy\'e sur un \'el\'ement de 
$E^+v_-\otimes v_-'+\St\otimes
(\St_{>-(p-1)\rho})$
et est donc non nul.
Il en d\'ecoule que tous les  $\tau(E^{(\nu)})v_-$
sont non nuls,
et qu'on doit donc avoir
$\tau(z_\nu)v_-=0$.
D'o\`u
$0=\chi_{-2(p-1)\rho}(\tau(z_\nu))
=\chi_{2(p-1)\rho}(z_\nu)$, comme voulu.

Supposons maintenant que
\begin{equation}
E^+F_i^{(rp)}=
F_i^{(rp)}
E^++
\underset{\substack{
s\in]0,rp[
\\
p\not\,| \ s}}{\sum}F_i^{(s)}x_s+
\sum_{\nu\in[0,p[^{R^+}}\sum_{s=0}^{r-1}F_i^{(sp)}z_{\nu s}y_{\nu s}
\end{equation}
pour certains
$x_s, y_{\nu s}\in \Dist(G)$
et
$z_{\nu s}\in\Dist(T_1)$
avec
$\chi_{2(p-1)\rho}(z_{\nu s})=0$.
Pour montrer qu'une formule similaire existe avec $r$ remplac\'e par  
$r+1$, il y a deux cas.
L'un d'entre eux est quand  
$n, m\leq r$ avec
$n+m=r+1$
tels que
$p\not|\binom{(r+1)p}{np}$.
Dans ce cas, des formules similaires \`a (3) existent avec $r$ remplac\'e par 
 $n$ et $m$.
Alors
\begin{align*}
\binom{(r+1)p}{np}
&
E^+F_i^{((r+1)p)}
=
E^+F_i^{(np)}F_i^{(mp)}
\\
&=
F_i^{(np)}
E^+F_i^{(mp)}
+
(\underset{\substack{
s\in]0,np[
\\
p\not\,| \ s}}{\sum}F_i^{(s)}x_s
+
\sum_{\nu\in[0,p[^{R^+}}\sum_{s=0}^{n-1}F_i^{(sp)}z_{\nu s}y_{\nu s})F_i^{(mp)}
\\
&=
F_i^{(np)}
(F_i^{(mp)}
E^+
+
\underset{\substack{
s\in]0,mp[
\\
p\not\,| \ s}}{\sum}F_i^{(s)}x_s'+
\sum_{\nu\in[0,p[^{R^+}}\sum_{s=0}^{m-1}F_i^{(sp)}z_{\nu s}'y_{\nu s}')
+
\\
&\hspace{3cm}
(\underset{\substack{
s\in]0,np[
\\
p\not\,| \ s}}{\sum}F_i^{(s)}x_s
+
\sum_{\nu\in[0,p[^{R^+}}\sum_{s=0}^{n-1}F_i^{(sp)}z_{\nu s}y_{\nu s})F_i^{(mp)}
\\
&=
\binom{(r+1)p}{np}
F_i^{((r+1)p)}
E^+
+
\underset{\substack{
s\in]0,(r+1)p[
\\
p\not\,| \ s}}{\sum}F_i^{(s)}x_s''
+
\sum_{\nu\in[0,p[^{R^+}}\sum_{s=0}^{r}F_i^{(sp)}z_{\nu s}''y_{\nu s}''
\end{align*}
pour certains
$x_s'', y_{\nu s}''\in \Dist(G)$
et
$z_{\nu s}''\in\Dist(T_1)$
avec
$\chi_{2(p-1)\rho}(z_{\nu s}'')=0$,
et (3) est v\'erifi\'e avec  $r$ remplac\'e par  
$r+1$.
Le second cas est quand il n'existe pas de tels  
$n$ et $m$.
Dans ce cas
$r+1$ est impair si  $p$ est impair,
et nous allons raisonner comme dans le premier cas lorsque  $r=1$.
Posons
$a=r+1$.
Alors
$\Dist(G_1)
\ni
\sum_{k=0}^{ap}
F_i^{(k)}E^+(-1)^{ap-k}F_i^{(ap-k)}
=
-E^+F_i^{(ap)}+F_i^{(ap)}E^+
+
\sum_{k=1}^{ap-1}
F_i^{(k)}E^+(-1)^{ap-k}F_i^{(ap-k)}$.
Si
$p|k$, soit
$k=bp$, on peut \'ecrire par  hypoth\`ese de r\'ecurrence
  \begin{align*}
F_i^{(bp)}E^+F_i^{(ap-bp)}
&=
F_i^{(bp)}
(F_i^{((a-b)p)}
E^++
\underset{\substack{
s\in]0,(a-b)p[
\\
p\not\,| \ s}}{\sum}F_i^{(s)}x_s
+
\sum_{\nu\in[0,p[^{R^+}}\sum_{s=0}^{a-b-1}F_i^{(sp)}z_{\nu s}y_{\nu s})
\end{align*}
pour certains
$x_s, y_{\nu s}\in \Dist(G)$
et
$z_{\nu s}\in\Dist(T_1)$
avec
$\chi_{2(p-1)\rho}(z_{\nu s})=0$.
Alors
\[
\Dist(G_1)
\ni
-E^+F_i^{(ap)}+F_i^{(ap)}E^+
+
\underset{\substack{
s\in]0,ap[
\\
p\not\,| \ s}}{\sum}F_i^{(s)}x_s
+
\sum_{\nu\in[0,p[^{R^+}}\sum_{s=1}^{a-b-1}F_i^{(sp)}z_{\nu s}y_{\nu s},
\]
et donc l'on peut \'ecrire comme dans (1)
 \begin{multline*}
E^+F_i^{(ap)}
=
F_i^{(ap)}E^+
+
\underset{\substack{
s\in]0,ap[
\\
p\not\,| \ s}}{\sum}F_i^{(s)}x_s
+
\sum_{\nu\in[0,p[^{R^+}}\sum_{s=1}^{a-b-1}F_i^{(sp)}z_{\nu s}y_{\nu s}
+
\sum_{j=1}^{p-1}F^{(j)}x_j'+\sum_{\nu\in[0,p[^{R^+}}
z_\nu E^{(\nu)}
\end{multline*}
pour un certain
$z_\nu\in\Dist(T_1)$.
Le m\^eme argument que pour   (1)
montre alors que
$\chi_{2(p-1)\rho}(z_\nu)=0$
pour tout $\nu$;
si
$F_i^{(ap)}\tau(E^+)v_-\ne0$,
$E_i^{(ap-1)}v_-\in\Delta(2(p-1)\rho)_{-2(p-1)\rho+(ap-1)\alpha_i}\setminus0$.
Alors
$s_i(-2(p-1)\rho+(ap-1)\alpha_i)=
-2(p-1)\rho+2(p-1)\alpha_i-(ap-1)\alpha_i=
-2(p-1)\rho-((a-2)p+1)\alpha_i$
serait aussi un poids de  $\Delta(2(p-1)\rho)$, ce qui est 
absurde.
Cela termine la v\'erification de l'induction.

\setcounter{equation}{0}
\noindent
(3.8)
{\bf Proposition}
(cf. \cite[Prop. 4.6, Th. 4.7]{KuLi}). 
 {\it {Pour tout   
$M\in\cC_B$,
d\'efinissons 
$\Psi_{2(p-1)\rho,M} : \{\cF_B^G(2(p-1)\rho\otimes
M^{[1]})\}^{\phi}\to
\cF_B^GM$
par 
$(\Psi_{2(p-1)\rho,M}f)(x)=f(E^+\phi(x))\otimes1$
pour tous $ f\in\{\cF_B^G(2(p-1)\rho\otimes
M^{[1]})\}^{\phi}$,
$x\in\Dist(G)$,
o\`u 1 dans le membre de droite est un \'el\'ement de base de   $-2(p-1)\rho$.
L'application est 
$\Psi_{2(p-1)\rho,M}$ est $G$-lin\'eaire,
i.e.,
$\Dist(G)$-lin\'eaire,
fonctorielle en
$M$,
et induit des applications
$G$-lin\'eaires naturelles
$\bfR^\bullet\Psi_{2(p-1)\rho,M} : \{\bfR^\bullet\cF_B^G(2(p-1)\rho\otimes
M^{[1]})\}^{\phi}\to
\bfR^\bullet\cF_B^GM$.

}}

\pf
Cela r\'esulte de 
(3.7) comme dans  \cite{KuLi}.

\setcounter{equation}{0}
\noindent
(3.9) Pour tous 
$ M_1, M_2\in\cC_P$,
l'application  $G$-lin\'eaire de
cup-produit
$\smallsmile : (\cF_P^GM_1)\otimes(\cF_P^GM_2)\to
\cF_P^G(M_1\otimes M_2)$
est d\'efinie par
$h_1\otimes h_2\mapsto(h_1\otimes h_2)\circ\Delta$
avec
$\Delta$ le coproduit sur $\Dist(G)$ et s'ins\`ere 
dans un diagramme commutatif
\[
\xymatrix
{
(\cF_P^GM_1)\otimes(\cF_P^GM_2)
\ar@{-->}[rr]^-{\smallsmile}
\ar@{-}[d]^\sim_-{\text{"identit\'e"}}
\ar
@{}
[]!<0ex,-7ex>;[r]
|{\hspace{0.5cm}\circlearrowleft}
&&
\cF_P^G(M_1\otimes M_2).
\\
\cF_P^G((\cF_P^GM_1)\otimes
M_2)
\ar!<0ex,-2.5ex>;[urr]_-{\quad\cF_P^G(\ev_{M_1}\otimes Id_{M_2})}
}
\]

\noindent
{\bf Lemme}
(cf. \cite[Lem. 4.9]{KuLi}).
{\it{Pour tous  $M_1, M_2\in\cC_P$,
le cup-produit induit des applications   $G$-lin\'eaires naturelles
$\smallsmile : (\cF_P^GM_1)\otimes(\bfR^\bullet\cF_P^GM_2)\to
\bfR^\bullet\cF_P^G(M_1\otimes M_2)$.}}

\setcounter{equation}{0}
\noindent
(3.10)
Soient
$v_-\in\Delta(2(p-1)\rho)_{-2(p-1)\rho}\setminus0$
un vecteur de plus bas poids, et
$f_0\in\nabla(2(p-1)\rho)_0$
avec
$f_0(E^+v_-)\ne0$
via la dualit\'e
$\nabla(2(p-1)\rho)\simeq\Delta(2(p-1)\rho)^*$.

\noindent
{\bf Th\'eor\`eme}
(cf. \cite[Prop. 4.11]{KuLi}). 
{\it
Consid\'erons
$f_0$
comme un \'el\'ement de
$\cF_B^G(2(p-1)\rho)$.
Alors, pour tout $M\in\cC_B$, on a un diagramme commutatif
 
\[
\xymatrix@R2ex{
((\bfR^\bullet\cF_B^GM)^{[1]})^\phi
\ar@{=}[ddd]
\ar[rr]^-{(\bfR^\bullet\Phi_M)^\phi}
&&
(\bfR^\bullet\cF_B^G(M^{[1]}))^\phi
\ar[dd]^-{(f_0\smallsmile?)^\phi}
\\
&&
\\
&&
(\bfR^\bullet\cF_B^G(2(p-1)\rho\otimes M^{[1]}))^\phi
\ar@{-}[d]^\sim
\\
\bfR^\bullet\cF_B^GM
&&
\bfR^\bullet\{(\cF_B^G(2(p-1)\rho\otimes M^{[1]}))^\phi\}
\ar[ll]^-{\bfR^\bullet\Psi_{2(p-1)\rho,M}},
}
\]
dans lequel
$f_0\smallsmile?$
est $T$-lin\'eaire mais non n\'ecessairement 
  $G$-lin\'eaire en g\'en\'eral.}

\begin{center}
{\bf
$4^\circ$
Faisceautisation}
\end{center}

Nous allons maintenant \'etendre les r\'esultats obtenus sur  $\bbF_p$ \`a un corps alg\'ebriquement clos   $\Bbbk$, puis faisceautiser les constructions pr\'ec\'edentes. Soient $G_\Bbbk=G\otimes_{\bbF_p}\Bbbk$ et $F_\Bbbk : G_\Bbbk\to G_\Bbbk$ le morphisme de Frobenius absolu de $G_\Bbbk$. On dispose  donc  de deux morphismes $F_\Bbbk^\sharp=Fr^\sharp\otimes_{\bbF_p}\Bbbk?^p : \Bbbk[G_\Bbbk]=\bbF_p[G]\otimes_{\bbF_p}\Bbbk\to \bbF_p[G]\otimes_{\bbF_p}\Bbbk$ et  $\Dist(F_\Bbbk)=\Dist(Fr)\otimes_{\bbF_p}?^{\frac{1}{p}}:\Dist(G_\Bbbk)=\Dist(G)\otimes_{\bbF_p}\Bbbk\to\Dist(G)\otimes_{\bbF_p}\Bbbk$ qui ne sont  pas $\Bbbk$-lin\'eaires. Posons $\phi_\Bbbk=\phi\otimes_{\bbF_p}?^p :\Dist(G_\Bbbk)=\Dist(G)\otimes_{\bbF_p}\Bbbk\to\Dist(G)\otimes_{\bbF_p}\Bbbk$.  Fixons d'autre part un sous-groupe parabolique  $P$ de $G$ et notons $P_\Bbbk=G\otimes_{\bbF_p}\Bbbk$.

\setcounter{equation}{0}
\noindent
(4.1)
Pour tout $\lambda\in\Lambda^+_P$, on pose 
$\nabla_{\Bbbk}(\lambda)=\nabla
(\lambda)\otimes_{\bbF_p}\Bbbk$, et l'on munit 
 $A
=\coprod_{\lambda\in\Lambda^+}\nabla_{\Bbbk}(\lambda)
$
  de la structure   de 
$G_\Bbbk
$-alg\`ebre  induite par le cup-produit.
On d\'efinit une application additive 
$\Psi_{A
}\in\Ab(A
,A
)$ par
\[
\xymatrix{
A
\ar@{-->}[rr]^-{\Psi_{A
}}
\ar@{}
;[dddrr]
|{\hspace{0cm}\circlearrowleft}
&&
A
\\
(\cF_P^G\lambda)\otimes_{\bbF_p}\Bbbk
\ar@{^(->}!<0ex,2.5ex>;[u]
\ar[d]
\\
\{\cF_P^G((\lambda^\phi)^{[1]})\}\otimes_{\bbF_p}\Bbbk
\ar[d]_{\mu_0\otimes_{\bbF_p}Id_{\Bbbk}
}
\\
\{(\cF_P^G((\lambda^\phi)^{[1]}))\otimes_{\bbF_p}\Bbbk\}
^{\phi_\Bbbk
}
\ar[rr]_-{\Psi_{\lambda}\otimes_{\bbF_p}?^{\frac{1}{p}}
}
&&
\{\cF_P^G(\lambda^\phi)
\}\otimes_{\bbF_p}\Bbbk
,
\ar@{^(->}!<0ex,2.5ex>;[uuu]
}
\]
avec
$(\cF_P^G\lambda)\otimes_{\bbF_p}\Bbbk
\to
\{\cF_P^G((\lambda^\phi)^{[1]})\}\otimes_{\bbF_p}\Bbbk
$
l'application identit\'e (resp. nulle)  si
$\lambda\in p\Lambda$
(resp. sinon).
L'application 
$\Psi_{A
}$ n'est pas $\Bbbk$-lin\'eaire, mais Frobenius-lin\'eaire au sens o\`u  pour tous 
$a, b\in A
$, on a
\[
\Psi_{A
}(a^pb)=a\Psi_{A
}(b),
\]
et scinde l'application d'\'el\'evation \`a la puissance  $p$ ; on a aussi
$\Psi_{A
}(1)=1$.
Notons
$\End_\Fr(A
)$ le $\Bbbk$-espace vectoriel des endomorphismes additifs Frobenius-lin\'eaires de $A
$. 
En  utilisant la structure de  $G%
_\Bbbk
$-alg\`ebre sur
$A
$,
on peut d\'efinir une action de   $G%
_\Bbbk
$ sur $\End_\Fr(A
)$
par
$g\bullet\psi=g\psi(g^{-1}?)$
pour tous $g\in G%
_\Bbbk
$,
$\psi\in\End_\Fr(A
)$
\cite[4.1.1.4]{BK}.
L'application 
$\Psi_{A
}$ est alors, par construction, $G%
_\Bbbk
$-semi-invariante (terminologie   sugg\'er\'ee par  \cite[4.15]{BK})
c'est-\`a-dire
$T%
_\Bbbk
$-lin\'eaire et telle que pour chaque sous-groupe radiciel  $\{x_{\pm\alpha_i}(\xi)\mid\xi\in\Bbbk\}$,
$i\in[1,\ell]$,
de
$G%
_\Bbbk
$, on ait 
\begin{align*}
x_{\alpha_i}(\xi)\bullet\Psi_{A
}=
\Psi_{A
}(\sum_{r\in[0,p[}(-\xi)^rE_i^{(r)}?)
\quad\text{et}\quad
x_{-\alpha_i}(\xi)\bullet\Psi_{A
}=
\Psi_{A
}(\sum_{r\in[0,p[}(-\xi)^rF_i^{(r)}?).
\end{align*}
La notion de
semi-invariance a \'et\'e originellement introduite par  Mathieu
\cite{M} (r\'ef\'erence dans laquelle une application 
  $B_\Bbbk$-semi-invariante, i.e. $T_\Bbbk$-lin\'eaire et v\'erifiant la formule ci-dessus   pour les sous-groupes radiciels
  $\{x_{-\alpha_i}(\xi)\mid\xi\in\Bbbk\}$,
$i\in[1,\ell]$,
de
$B%
_\Bbbk
$, est dite 
$B_\Bbbk$-canonique).
Pour une
$T%
_\Bbbk
$-alg\`ebre $Q$ et pour
$\psi\in\End_\Fr(Q)$,
l'application $\psi$ est
$T%
_\Bbbk
$-lin\'eaire si et seulement si  pour tout
$\lambda\in\Lambda$,
\[
\psi(Q_{\lambda})\subseteq
\begin{cases}
Q_{\lambda^1}
&\text{si $\lambda\in p\Lambda$}
\\
0&\text{sinon}.
\end{cases}
\]

Notons que pour tous 
$g\in G%
_\Bbbk
$, $a\in A
$, on a 
$g(a^p)=(ga)^p$,
trivialisant ainsi l'action du noyau de Frobenius   $G_{\Bbbk,1}$;
si
$a_1,\dots, a_n$ est une base sur  $\Bbbk$ 
d'un sous-espace  $G%
_\Bbbk
$-invariant de
$A
$,
alors les
$a_1^p,\dots, a_n^p$
forment une base d'un autre sous-espace    $G%
_\Bbbk
$-invariant avec
$g(a_i^p)=(\sum_jg_{ji}a_j)^p
=\sum_jg_{ji}^pa_j^p$,
et donc la 
  $G%
_\Bbbk
$-action sur 
$\sum_i%
\Bbbk 
a_i^p$
se factorise \`a  travers le   Frobenius.
En particulier, les  actions des
$E_i^{(r)}$ et des  $F_i^{(r)}$,
$i\in[1,\ell]$, $r\in[0,p[$, sur
$A
$ sont toutes Frobenius lin\'eaires.
Ainsi

\noindent
{\bf Th\'eor\`eme}
(cf. \cite[Prop. 6.2]{KuLi}). 
{\it{L'application $\Psi_{A
}$
est un scindage Frobenius lin\'eaire  $G%
_\Bbbk
$-semi-invariant de l'application d'\'el\'evation \`a la puissance $p$ sur 
$A
$}}.

\setcounter{equation}{0}
\noindent
(4.2) On va maintenaint faisceautiser l'application $\Psi_{A
}$.
Posons
$\cP=G%
_\Bbbk
/P%
_\Bbbk
$ et $\Lambda_P^{++}=\{\lambda \in \Lambda_P\mid \langle\lambda,
\alpha^\vee\rangle>0 \ {\text {pour tout}}\ \alpha\in R^s\setminus R_P^s\}$.
On dispose   [II.8.5]{J}  d'une immersion ferm\'ee $\cP\to
\bbP(\Delta_\Bbbk(-w_0\lambda))$
($\Delta_\Bbbk(-w_0\lambda)=\nabla_\Bbbk(\lambda)^*$)
via
$g\mapsto[gv_-]$,
avec
$v_-$
un vecteur de plus bas poids de 
$\Delta_\Bbbk(-w_0\lambda)$. D\'efinissons une 
 $\Bbbk$-alg\`ebre (en utilisant le cup-produit) gradu\'ee par
 
$\Gamma^\bullet(\cP,\lambda)=$
$\coprod_{m\in\bbN}\Gamma(\cP,i^*\cO_{\bbP(\Delta_\Bbbk(-w_0\lambda))}(m))\simeq$
$\coprod_{m\in\bbN}\Gamma(\cP,\cL
(m\lambda))
=
\coprod_{m\in\bbN}\nabla_{\Bbbk}(m\lambda)$.
Alors
$\cP\simeq
\Proj(\Gamma^\bullet(\cP,\lambda))$
par
\cite[Ex. II.5.14b+Ex. II.2.14c]{H}.
Remarquons que 
$\Gamma^\bullet(\cP,\lambda)$
est int\`egre car
$\nabla%
_{\Bbbk}
(m\lambda)=\Sch_{\Bbbk}(G%
_\Bbbk
,m\lambda)^{P%
_\Bbbk
}$.
Pour tout $w\in W$,
choisissons
$f_w\in\nabla
(\lambda)_{w\lambda}\setminus0$,
et posons
$\cP_w=\{x\in\cP\mid
f_w(x)\ne0\}$.
Alors
$\cP_w=wU%
_{P,\Bbbk}
^+P%
_\Bbbk
/P%
_\Bbbk
$
est un ouvert affine de 
$\cP$
avec
$\Bbbk[\cP_w]
\simeq
\Gamma^\bullet(\cP,\lambda)_{(f_w)}
=
\displaystyle\cup_{m\in\bbN}
\{
\frac{f}{{f_w}^m}\mid
f\in\nabla%
_\Bbbk
(m\lambda)\}$
dans
$\Frac(\Gamma^\bullet(\cP,\lambda))$.
En particulier,
les ouverts principaux d\'efinis par les  $f_w$, $w\in W$,
forment un recouvrement de
$\cP$. 
Suivant
\cite[II.5.15]{H},
le faisceau $F_*\cO_\cP$ se r\'ecup\`ere \`a partir   du 
 $\Gamma^\bullet(\cP, \lambda)$-module gradu\'e
$\Gamma^\bullet(\cP, F_*\cO_\cP)=
\coprod_{m\in\bbN}\Gamma(\cP,
(F_*\cO_\cP)(m))\simeq
\coprod_{m\in\bbN}\Gamma(\cP,
\cL(pm\lambda))$.
On a
$\Gamma^\bullet(\cP,
F_*\cO_\cP)_{(f_w)}
=
\cup_{m\in\bbN}
\{
\frac{f}{{f_w}^{pm}}\mid
f\in\nabla%
_\Bbbk
(pm\lambda)\}$.
Posons pour all\'eger
$A_w(\lambda)=\Gamma^\bullet(\cP,
\lambda)_{(f_w)}$
et
$A'_w(\lambda)=\Gamma^\bullet(\cP,
F_*\cO_\cP)_{(f_w)}$.
Definissons maintenant $\Psi_w^\lambda : A'_w(\lambda)\to
A_w(\lambda)$
par
$\displaystyle
\frac{f}{{f_w}^{pm}}\mapsto
\frac{\Psi_{A
}(f)}{{f_w}^{m}}$,
$f\in\nabla_\Bbbk(pm\lambda)$.
C'est bien d\'efini et les 
  $\Psi_w^\lambda$,
$w\in W$, se recollent ensemble pour donner un morphisme  $\Theta : F_{
*}\cO_\cP\to
\cO_\cP$,
qui est en fait ind\'ependant du choix de  $\lambda\in\Lambda_P^{++}$.
Par construction

\noindent
{\bf Th\'eor\`eme}
(cf. \cite[Th. 9.6]{KuLi}). 
{\it L'application
$\Theta : F_{
*}\cO_\cP\to
\cO_\cP$ est un scindage  $G%
_\Bbbk
$-semi-invariant
du comorphisme
$F
^\sharp :
\cO_\cP\to
F_{
*}\cO_\cP$
d\'efinissant l'endomorphisme de   Frobenius absolu de $\cP$.

}

\setcounter{equation}{0}
\noindent
(4.3)
On peut g\'en\'eraliser le th\'eor\`eme pr\'ec\'edent \`a la situation suivante. Soient
$X$ un $\Bbbk$-sch\'ema muni de son endomorphisme de 
Frobenius absolu
$F_X$ et
$\cL$ un faisceau inversible sur  $X$ et
$f\in\Gamma(X,\cL)\simeq\Mod_X(\cO_X,\cL)$.
On dira que $X$  est Frobenius
$f$-scind\'e si     le compos\'e
$\cO_X\overset{F_X^\sharp}{\longrightarrow}
F_{X*}\cO_X
\overset{F_{X*}f}{\longrightarrow}
F_{X*}\cL$
admet un inverse \`a gauche.
On suppose maintenant que  $P=B$
et l'on \'ecrira 
$\cB$ \`a la place de  $\cP$.
Utilisant cette fois  (3.10)
on obtient de m\^eme

\noindent
{\bf Th\'eor\`eme}
(cf. \cite[Th. 6.5]{KuLi}). 
{\it
Soit
$f_0\in\nabla
(2(p-1)\rho)_0\setminus0$
vu comme morphisme de $\cO_\cB$
vers
$\cL_\cB(2(p-1)\rho)$.
Alors le compos\'e
$F_{
*}(f_0)\circ
F
^\sharp : \cO_\cB\to F_{
*}\cL_\cB(2(p-1)\rho)$
admet un inverse \`a gauche $\cO_\cB$-lin\'eaire et $G_\Bbbk$-semi-invariant.}

\pf
Soit
$\lambda\in\Lambda_B^{++}$
comme dans (4.2).
Posons
$A'(2(p-1)\rho,\lambda)=
\Gamma^\bullet(\cB,F_*\cL(2(p-1)\rho))$,
qui redonne
$F_*\cL(2(p-1)\rho)$
par faisceautisation.
Definissons
$\Psi_{2(p-1)\rho}^\lambda :
A'(2(p-1)\rho,\lambda)\to
A(\lambda)$
via le diagramme
\[
\xymatrix
@R2ex
{
A'(2(p-1)\rho,\lambda)
\ar@{-->}[rrrr]
\ar
@{}
[]!<4ex,0ex>;[dddrrrr]
|{\circlearrowleft}
&&&&
A(\lambda)
\\\\
\cF_B^G(2(p-1)\rho+pm\lambda)\otimes_{\bbF_p}\Bbbk
\ar@{^(->}[]!<0ex,3ex>;[uu]
\ar@{=}[d]
&&&&
\\
\cF_B^G(2(p-1)\rho+(m\lambda)^{[1]})^\phi
\otimes_{\bbF_p}\Bbbk
\ar[rrrr]_-{\Psi_{2(p-1)\rho,m\lambda}\otimes_{\bbF_p}?^{\frac{1}{p}}
}
&&&&
\cF_B^G(m\lambda)\otimes_{\bbF_p}\Bbbk
\ar@{^(->}[]!<0ex,3ex>;[uuu],
}
\]
  dans lequel l'identification
$\cF_B^G(2(p-1)\rho+pm\lambda)=\cF_B^G(2(p-1)\rho+(m\lambda)^{[1]})^\phi
$
est faite en tant que $\bbF_p$-espaces vectoriels. Pour tous
$a\in\nabla_\Bbbk(n\lambda)$,
$b\in\nabla_\Bbbk(2(p-1)\rho+pm\lambda)$,
on a de nouveau
$\Psi_{2(p-1)\rho,(n+m)\lambda}
(a^pb)=a\Psi_{2(p-1)\rho,m\lambda}
(b)$.
Pour tout $w\in W$,
posons
$A'_w(2(p-1)\rho,\lambda)=
\Gamma^\bullet(\cB,F_*\cL(2(p-1)\rho))_{(f_w)}=
\cup_{m\in\bbN}
\{
\frac{f}{{f_w}^{pm}}\mid
f\in\nabla%
_\Bbbk
(2(p-1)\rho+pm\lambda)\}$,
et d\'efinissons
$\Psi^\lambda_{2(p-1)\rho,w} :
A'_w(2(p-1)\rho,\lambda)\to
A_w(\lambda)$
via
$\frac{f}{{f_w}^{pm}}\mapsto
\frac{\Psi_{2(p-1)\rho}^\lambda(f)}{{f_w}^{m}}
$
pour
$f\in\nabla_\Bbbk(2(p-1)\rho+pm\lambda)$.
Les
$\Psi_{2(p-1)\rho,w}^\lambda$, $w\in W$, sont bien d\'efinis et se recollent ensemble pour former un morphisme $\cO_\cB$-lin\'eaire 
$F_*\cL(2(p-1)\rho)\to \cO_\cB$
qui scinde
$F_*(f_0)\circ F^\sharp$.

La $G_\Bbbk$-semi-invariance de
$\Psi_{2(p-1)\rho}^\lambda$ r\'esulte alors de la  Frobenius-lin\'earit\'e et de la    $G$-lin\'earit\'e de
$\Psi_{2(p-1)\rho,m\lambda}$.

\setcounter{equation}{0}
\noindent
(4.4) On peut \'egalement chercher \`a scinder des sous-sch\'emas :
pour un   $\Bbbk$-sous-sch\'ema $Y$ d'un
  $\Bbbk$-sch\'ema $X$
d\'efini par un faisceau d'id\'eaux
$\cI_Y$
on dira que que  $Y$ est Frobenius-scind\'e de mani\`ere compatible si et seulement si il existe un scindage de Frobenius $\sigma$ de $X$ tel que  
$\sigma(F_*\cI_Y)\subseteq\cI_Y$
\cite[1.1.3]{BK}, auquel cas on dira simplement que $\sigma$ scinde de mani\`ere compatible $Y$.

Pour tout $w\in W$, soient
$X
(w)=\overline{B_\Bbbk wB_\Bbbk/P_\Bbbk}$
et
$X^+
(w)=\overline{B^+_\Bbbk
wB_\Bbbk/P_\Bbbk}$,
les sous-sch\'emas de   Schubert de $\cP$.
On va montrer que  tous ceux-ci sont scind\'es de mani\`ere compatible par  l'application
$\Theta$ de (4.2).

Comme dans  (4.2) prenons
$\lambda\in\Lambda^{+}_P$
avec
$\langle\lambda,\alpha^\vee\rangle>0$
pour tout $\alpha\in R^s\setminus R^s_P$,
et choisissons   
$v_-\in\Delta(-w_0m\lambda)_{-m\lambda}\setminus0$ pour tout $m\in\bbN^+$.
Rappelons 
\cite[II.14.19]{J} qu'il existe un
isomorphisme de $G$-modules
\begin{equation}
\nabla(m\lambda)\simeq
\Delta(-mw_0\lambda)^*
\simeq
\cF_P^Gm\lambda
\quad\text{via}\quad
h(?v_-)\longleftarrow\!\shortmid
h
\longmapsto
h(S(?)v_-),
\end{equation}
avec $S$ l'antipode de
$\Dist(G)$.
Soit
$\cI(w)$ le faisceau d'id\'eaux de
$\cO_\cP$ d\'efinissant $X(w)$.
Si $i_w : X
(w)\hookrightarrow \cP$ d\'esigne l'immersion canonique,
on a une suite exacte  $0\to\cI(w)\to\cO_\cP\to
i_{w*}\cO_{X(w)}\to0$.
Comme dans (4.3) on a
$\cI(w)=\Gamma^\bullet(\cP, \cI(w))^\sim$
et
$i_{w*}\cO_{X(w)}=\Gamma^\bullet(\cP, i_{w*}\cO_{X(w)})^\sim$
avec
$\Gamma^\bullet(\cP, \cI(w))
=
\coprod_{m\in\bbN}\Gamma(\cP,\cI(w)(m))
$.
Comme
$
\Gamma^\bullet(\cP, i_{w*}\cO_{X(w)})
\simeq
\coprod_{m\in\bbN}\Sch_{\Bbbk}(\overline{B_\Bbbk
wP_\Bbbk},
m\lambda)^{P_\Bbbk}$,
on a une suite exacte
\[
0\to
\Gamma(\cP, \cI(w)(m))\to
\nabla_\Bbbk(m\lambda)\overset{\res}{\longrightarrow}
\Sch_{\Bbbk}(\overline{B_\Bbbk
wP_\Bbbk},
m\lambda)^{P_\Bbbk}\to0.
\]
Posons
$\nabla_w(m\lambda)=\Gamma(X(w),i_w^*\cL(m\lambda))=
\Sch_{\Bbbk}(\overline{B_\Bbbk
wP_\Bbbk},
m\lambda)^{P_\Bbbk}$.
D'un autre c\^ot\'e,
si $v_-\in\Delta(-mw_0\lambda)_{-m\lambda}\setminus0$
et si
$(\Dist(U)wv_-)^\perp=\{
h\in\Delta_\Bbbk(-mw_0\lambda)^*\mid
h(\Dist(U)wv_-)=0\}$,
on tire de
\cite[II.14.19.2, 3]{J} 
une suite exacte
$0\to
(\Dist(U)wv_-)^\perp\to
\nabla_\Bbbk(m\lambda)
\overset{\res}{\longrightarrow}
\nabla_w(m\lambda)$,
et donc
\begin{equation}
\Gamma(\cP, \cI(w)(m))\simeq(\Dist(U)wv_-)^\perp.
\end{equation}
De m\^eme avec
$i_w^+ : X^+(w)\hookrightarrow 
\cP$,
 si l'on pose
$\nabla_w^+(m\lambda)=\Gamma(X^+(w), i_w^*\cL(m\lambda))$, on a   une suite exacte  
\[
0\to
(\Dist(U^+)wv_-)^\perp\to
\nabla_\Bbbk(m\lambda)
\overset{\res}{\longrightarrow}
\nabla_w^+(m\lambda).
\]

\noindent
{\bf Th\'eor\`eme}
(cf. \cite[Th. 6.7]{KuLi}). 
{\it
Soit
$Y$ un sous-$\Bbbk$-sch\'ema de
$\cP$ obtenu  \`a partir des 
$X
(w)$
et des
$X
^+(w)$,
$w\in W$, en combinant   des unions ou intersections sch\'ematiques et des consid\'erations de composantes irr\'eductibles r\'eduites.
Alors $Y$ est scind\'e de mani\`ere compatible par  $\Theta$.
En particulier, $Y$ est r\'eduit.

}

\pf
On va d'abord montrer que
$X
(w)$ est scind\'e de mani\`ere compatible par
$\Theta$ et ce pour tout
$w\in W$.
Remarquons tout d'abord que 
$\cP_w\cap X(w)\ni wP_\Bbbk$.
Si
$v_{-m\lambda}\in\Delta(-mw_0\lambda)_{-m\lambda}\setminus0$,
on a
\begin{align*}
\Gamma^\bullet(\cP,
F_*\cI(w))
&=
\coprod_{m\in\bbN}
\Gamma(\cP,(F_*\cI(w))(m))
\simeq
\coprod_{m\in\bbN}
\Gamma(\cP,\cI(w)\otimes_{O_{\cP}}
\cL(pm\lambda))
\\
&\simeq
\coprod_{m\in\bbN}
(\Dist(U)wv_{-pm\lambda})^\perp
\quad\text{par (2)},
\end{align*}
et donc
$\Gamma(\cP_w,F_*\cI(w))=
\Gamma^\bullet(\cP,
F_*\cI(w))_{(f_w)}=\cup_{m\in\bbN}\{\frac{f}{{f_w}^{pm}}
\mid
\tilde f\in(\Dist(U^+)wv_{-pm\lambda})^\perp\}$,
avec
$\tilde f\in\Delta_\Bbbk(-pmw_0\lambda)^*$ correspondant \`a
$f\in\nabla_\Bbbk(pm\lambda)$ via (1).
Comme le morphisme 
$\Gamma(\cP_w,\Theta) :
\Gamma(\cP_w, F_*\cO_\cP)=
\Gamma^\bullet(\cP_w, F_*\cO_\cP)_{(f_w)}
=
\cup_{m\in\bbN}\{
\frac{f}{{f_w}^{pm}}
\mid
f\in\nabla_\Bbbk(pm\lambda)
\}\to
\Bbbk[\cP_w]=
\cup_{m\in\bbN}\{
\frac{f}{{f_w}^{m}}
\mid
f\in\nabla_\Bbbk(m\lambda)
\}$
se d\'ecrit comme
$\frac{f}{{f_w}^{pm}}\mapsto
\frac{\Psi_A(f)}{{f_w}^m}$,
on a seulement \`a montrer que  
$\widetilde{\Psi_A(f)}=\widetilde{\Psi_{m\lambda}(\mu_0f)}\in
(\Dist(U)wv_{-m\lambda})^\perp$
pour tout $\tilde f\in(\Dist(U)wv_{-pm\lambda})^\perp$.

Si
$\tilde f\in(\Dist(U)zv_-)^\perp$,
alors pour tous
$t\in T$, $\mu\in\Dist(U)$,
\begin{align*}
(t\tilde f)(\mu zv_-)
&=
\tilde f(t^{-1}\mu zv_-)
=
\tilde
f(\Ad(t^{-1})(\mu)t^{-1}zv_-)
\\
&=0
\quad\text{car
$\Ad(t^{-1})(\mu)\in\Dist(U)$}.
\end{align*}
Donc
$(\Dist(U)zv_-)^\perp$ admet une d\'ecomposition suivant les poids et par suite on peut supposer que   
$f$ est vecteur propre associ\'e \`a un certain poids dans $p\Lambda$.
Posons $\tilde f'=\widetilde{\Psi_{m\lambda}(\mu_0f)}
$.

Pour tout $\mu\in\Dist(U)$, on a  donc
\begin{align*}
\tilde f'(\mu wv_{-m\lambda})
&=
\tilde f'(w(^{w^{-1}}\mu)
v_{-m\lambda})
\quad\text{avec}\quad
^{w^{-1}}\mu=\Ad(w^{-1})(\mu)
\\
&=
(w^{-1}\tilde f')((^{w^{-1}}\mu)
v_{-m\lambda})
=
(w^{-1}\tilde f')(SS(^{w^{-1}}\mu)
v_{-m\lambda}))
\\
&=
(w^{-1}\Psi_{m\lambda}(f))(S(^{w^{-1}}\mu)
)
\\
&=
\Psi_{m\lambda}(w^{-1}\bullet
f)(S(^{w^{-1}}\mu)
)
\quad\text{car
$\Psi_{m\lambda}$ est $G$-lin\'eaire par (3.4)}
\\
&=
\Psi_{m\lambda}(\phi(w^{-1})
f)(S(^{w^{-1}}\mu)
)
=
(\phi(w^{-1})f)(\phi(
S(^{w^{-1}}\mu)
))
\\
&=
(\phi(w^{-1})f)(S(
\phi(^{w^{-1}}\mu)
))
=
\widetilde{\phi(w^{-1})f}(\phi(^{w^{-1}}\mu)
v_-^{pm})
=
(\phi(w^{-1})\tilde f)(\phi(^{w^{-1}}\mu)
v_-^{pm})
\\
&=
\tilde f(\phi(w)\phi(^{w^{-1}}\mu)
v_-^{pm})
=
\tilde f(\phi(w^{w^{-1}}\mu)
v_-^{pm})
=
\tilde f(\phi(\mu w)
v_-^{pm})
=
\tilde f(\phi(\mu)\phi(w)
v_-^{pm})
\\
&=
\tilde f(\phi(\mu)w
v_-^{pm})
\quad\text{\`a  $\bbF_p^\times$ pr\`es
par (2.2) car
$\Delta(-w_0pm\lambda)_{-wpm\lambda}$ est de dimension 1}
\\
&=0
\quad\text{car
$\phi(\mu)\in\Dist(U)\mu_0$},
\end{align*}
comme voulu.

Comme $X(w)$ est int\`egre \cite[II.13.8]{J},
il d\'ecoule de
\cite[F.12]{J}
que
$\Theta(F_*\cI(w))\subseteq
\cI(w)$, et donc que
$X
(w)$ est scind\'e de mani\`ere compatible par
$\Theta$. De m\^eme $X
^+(w)$ est scind\'e de mani\`ere compatible par
$\Theta$.
Ainsi toute it\'eration d'une union ou d'une intersection sch\'ematique des   $X
(w)$
et des
$X
^+(w)$ est scind\'ee de mani\`ere compatible par $\Theta$, et de m\^eme pour les composantes irr\'eductibles r\'eduites \cite[F.13]{J}.

\setcounter{equation}{0}
\noindent
(4.5)
Montrons finalement  comment le scindage de 
\cite[Thm.(3.4)]{K95} s'interpr\`ete avec le point de vue que nous venons de d\'evelopper.
Supposons que $P=B$,  posons  $\cB=G_\Bbbk/B_\Bbbk$  et soit
$v_-\in\St_{-(p-1)\rho}\setminus0$
avec
$v_-^*\in(\St^*)_{(p-1)\rho}$
tel que
$v_-^*(v_-)=1$.
Alors le   $v_-$-scindage de Frobenius  de $\cB$ introduit dans  \cite{K95} et que nous noterons ici 
$\sigma: F_*\cL((p-1)\rho)\to\cO_\cB$
s'ins\`ere dans un diagramme   commutatif
\[
\xymatrix{
\cO_\cB
\ar[rr]^-{F^\sharp}
\ar[rd]
^-{v_-\otimes_{\bbF_p}Id_{\cO_\cB}}
&&
F_*\cO_\cB
\ar[dd]^-{F_*v_-}
\\
&
\qquad\St\otimes_{\bbF_p}\cO_\cB
\ar[dl]^-{v_-^*\otimes_{\bbF_p}Id_{\cO_\cB}}
\\
\cO_\cB
\ar@{=}[uu]
&&
F_*\cL((p-1)\rho).
\ar@{-}[ul]!<0ex,2ex>_-\sim
\ar[ll]^-{\sigma}
}
\]
De \cite{K95} (resp.
\cite[F.22]{J}) on tire que
$\sigma$ (resp. $\sigma\circ F_*v_-$)
scinde de mani\`ere compatible tous les sch\'emas de  Schubert  
$X^+(w)$ (resp. $X(w)$), $w\in W$, et donc
$\sigma\circ F_*v_-$ scinde de mani\`ere compatible tous les sous-sch\'emas   $Y$ comme dans  (4.4).

Donnons une v\'erification directe de ces r\'esultats;
\cite{K95}
(resp. \cite{J})
montre que le sous-sch\'ema ferm\'e  $\cB \setminus (U^+_\Bbbk
B_\Bbbk/B_\Bbbk)$
(resp. $\cB\setminus (U_\Bbbk
w_0
B_\Bbbk/B_\Bbbk)$)
est scind\'e de mani\`ere compatible. Prenons  $\lambda\in\Lambda^{+}$
avec
$\langle\lambda,\alpha^\vee\rangle>0$
pour tout  $\alpha\in R^s$
et posons
$A
=\coprod_{m\in\bbN}\nabla_\Bbbk(m\lambda)$,
$A'
=
(p-1)\rho\otimes_\Bbbk\coprod_{m\in\bbN}\nabla_\Bbbk((p-1)\rho+pm\lambda)$.
D\'efinissons une application additive 
$\Psi_{(p-1)\rho}^\lambda : A' \rightarrow A$
par la commutativit\'e du diagramme
\[
\xymatrix{
A'
\ar@{-->}[rrrrrr]^-{\Psi_{(p-1)\rho}^\lambda}
&&&&&&
A
\\
(p-1)\rho\otimes_\Bbbk\nabla_\Bbbk((p-1)\rho+pm\lambda)
\ar@{^(->}!<0ex,2.5ex>;[u]
\ar@{-}[d]_-\sim
&&&&&&
\\
(p-1)\rho\otimes_{\bbF_p}\St\otimes_{\bbF_p}\nabla(m\lambda)^{[1]}\otimes_{\bbF_p}\Bbbk
\ar[d]_-{\mu_0\otimes_{\bbF_p}\Bbbk}
\\
\{(p-1)\rho\otimes_{\bbF_p}\St\otimes_{\bbF_p}\nabla(m\lambda)^{[1]}\}^\phi\otimes_{\bbF_p}\Bbbk
\ar[rrrrrr]_-{\{(p-1)\rho\otimes_{\bbF_p}v_-^*\otimes_{\bbF_p}\nabla(m\lambda)^{[1]}\}^\phi\otimes_{\bbF_p}?^{\frac{1}{p}}}
&&&&&&
\nabla_\Bbbk(m\lambda).
\ar@{^(->}!<0ex,2.5ex>;[uuu]
}
\]
Comme $v_-^*:\St\to-(p-1)\rho$
est
$B_\Bbbk^+$-lin\'eaire,
c'est aussi le cas pour
$\Psi^\lambda_{(p-1)\rho}$.
D'autre part,
si $v\in\St_\eta$ avec
$\eta\ne-(p-1)\rho$, alors pour tous
$r\in\bbN^+$,
$z\in\nabla(m\lambda)^{[1]}$,
\[
F_i^{(pr)}(v\otimes z)\in\ker(v_-^*\otimes\nabla(m\lambda)^{[1]}).
\]
On a
$F_i^{(pr)}(v\otimes z)=
\sum_{k=0}^{pr}(F_i^{(pr-k)}v)\otimes
(F_i^{(k)}z)=
\sum_{k=0}^{r}(F_i^{(pr-pk)}v)\otimes
(F_i^{(k)}z)$.
Si $F_i^{(pk)}v\in\bbF_pv_-\setminus0$
il existe $k>0$ tel que 
$-(p-1)\rho+kp\alpha_i$
soit un poids de 
$\St$, et donc
$s_i(-(p-1)\rho+kp\alpha_i)=
-(p-1)\rho-(kp-(p-1))\alpha_i
<-(p-1)\rho$,
ce qui est absurde.
Il en d\'ecoule que pour tous  $v\in\St_\eta$,
$\eta\in\Lambda$,
on a 
\begin{align*}
(v_-^*\otimes\nabla(m\lambda)^{[1]})
&(F_i^{(pr)}(v\otimes z))
\\
&=
\begin{cases}
v_-^*(v)F_i^{(r)}z
=
F_i^{(r)}(v_-^*\otimes\nabla(m\lambda)^{[1]})(v\otimes z)
&\text{si
$\eta=-(p-1)\rho$},
\\
0
&\text{sinon},
\end{cases}
\end{align*}
et par suite,  pour tout  $\xi\in\Bbbk$,
\begin{align*}
x_{-\alpha_i}(\xi)
&\bullet\Psi^\lambda_{(p-1)\rho}(1\otimes
v\otimes z\otimes1)
=
x_{-\alpha_i}(\xi)
\Psi^\lambda_{(p-1)\rho}(x_{-\alpha_i}(-\xi)(1\otimes
v\otimes z\otimes1))
\\
&=
x_{-\alpha_i}(\xi)
\Psi^\lambda_{(p-1)\rho}(\sum_{r\in\bbN}(-\xi)^rF_i^{(r)}(1\otimes
v\otimes z\otimes1))
\\
&=
x_{-\alpha_i}(\xi)
\Psi^\lambda_{(p-1)\rho}(1\otimes
\sum_{r\in\bbN}F_i^{(r)}(v\otimes 
z)\otimes(-\xi)^{r})
\\
&=
x_{-\alpha_i}(\xi)
\Psi^\lambda_{(p-1)\rho}(1\otimes
\underset{\substack{r_0\in[0,p[
\\
r_1\in\bbN}}{\sum}
F_i^{(pr_1)}F_i^{(r_0)}(v\otimes 
z)\otimes(-\xi)^{r_0+pr_1})
\\
&=
x_{-\alpha_i}(\xi)
\sum_{r_1\in\bbN}(-\xi)^{r_1}F_i^{(r_1)}
\Psi^\lambda_{(p-1)\rho}(1\otimes
\sum_{r_0\in[0,p[}F_i^{(r_0)}(v\otimes
z)\otimes(-\xi)^{r_0})
\\
&=
\Psi^\lambda_{(p-1)\rho}(1\otimes
\sum_{r_0\in[0,p[}F_i^{(r_0)}(v\otimes
z)\otimes(-\xi)^{r_0}).
\end{align*}
Ainsi
$\Psi^\lambda_{(p-1)\rho}$
est
$B_\Bbbk^+$-lin\'eaire et $B_\Bbbk$-semi-invariant, et donc, en particulier,
$G_\Bbbk$-semi-invariant.
On remarquera, cependant, qu'aucun
 $(p-1)\rho\otimes_\Bbbk\nabla_\Bbbk((p-1)\rho+pm\lambda)$
n'est   muni d'une structure de  
$G_\Bbbk$-module;
on a d\'ecal\'e la structure de  $T_\Bbbk$-module sur le 
$G_\Bbbk$-module
$\nabla_\Bbbk((p-1)\rho+pm\lambda)$
par $(p-1)\rho$
pour rendre $T_\Bbbk$-lin\'eaire l'application 
$\Psi_{(p-1)\rho}^\lambda$.

Pour tout
$w\in W$, soit
$f_w\in\nabla(\lambda)_{w\lambda}\setminus0$
comme dans  (4.2) et \'ecrivons ici
$\cB_w$ pour $\cP_w$.
Si l'on d\'efinit
$\Psi_{(p-1)\rho,w}^\lambda:\cup_{m\in\bbN}\{\frac{f}{{f_w}^{pm}}\mid
f\in(p-1)\rho\otimes_\Bbbk\nabla_\Bbbk((p-1)\rho+pm\lambda)\}\to
\Bbbk[\cB_w]$
via
$\frac{f}{{f_w}^{pm}}\mapsto
\frac{\Psi_{(p-1)\rho}^\lambda(f)}{{f_w}^m}$,
alors les $\Psi_{(p-1)\rho,w}^\lambda$,
$w\in W$, se recollent pour former  $\sigma$.

Pour voir directement que les  
$\sigma$ scindent de mani\`ere compatible les  
$X^+(w)$, $w\in W$, i.e., que
$\sigma\{F_*(\cI^+(w)\otimes_\cB\cL((p-1)\rho)\}\subseteq
\cI^+(w)$ pour le faisceau d'id\'eaux 
$\cI^+(w)$ de
$X^+(w)$,
soit
$v_{-pm\lambda}^{(p-1)\rho}$ 
(resp.
$v_{-m\lambda}$) un vecteur de plus bas poids de   $\Delta(-w_0((p-1)\rho+pm\lambda))$
(resp.
$\Delta(-w_0m\lambda)$).
On doit voir l'existence d'un diagramme commutatif
\[
\xymatrix@R2ex{
\{(p-1)\rho\otimes\St\otimes
\nabla(m\lambda)^{[1]}\}^\phi
\ar[rrrr]^-{\{(p-1)\rho\otimes
v_-^*\otimes
\nabla(m\lambda)^{[1]}\}^\phi}
&&&&
\nabla(m\lambda)
\\\\
(p-1)\rho\otimes\St\otimes
\nabla(m\lambda)^{[1]}
\ar[uu]^-{\mu_0}
&&&&
(\Dist(U^+)wv_{-m\lambda})^\perp.
\ar@{^(->}!<0ex,2.5ex>;[uu]
\\
(p-1)\rho\otimes
\nabla((p-1)\rho+pm\lambda)
\ar@{-}[u]^-\sim
&&&&
\\\\
(p-1)\rho\otimes
(\Dist(U^+)wv_{-m\lambda}^{(p-1)\rho})^\perp
\ar@{^(->}!<0ex,2.5ex>;[uu]
\ar@{-->}!<0ex,-4ex>;[uuurrrr]
}\]
Rappelons tout d'abord l'isomorphisme   $G$-lin\'eaire  
$\nabla((p-1)\rho+pm\lambda)
\simeq
\St\otimes\nabla(m\lambda)^{[1]}$
via
\[
\sum v_i\smallsmile f_i^p\longmapsto
\sum v_i\otimes f_i,
\quad
v_i\in\St, f_i\in\nabla(m\lambda).
\]
Notons par  
$\tilde f\in
\Delta(-w_0((p-1)\rho+pm\lambda))^*$
l'\'el\'ement correspondant \`a   $f\in\nabla((p-1)\rho+pm\lambda)$.
On est ramen\'e \`a v\'erifier que
$\sum
v_-^*(v_i)\tilde f_i(\Dist(U^+)wv_{-m\lambda})=0$
si
$\widetilde{\sum v_i\smallsmile f_i^p}(
\Dist(U^+)wv_{-pm\lambda}^{(p-1)\rho})=0$,
i.e., 
$\sum
v_-^*(v_i)f_i(gw)=0$
pour tout $ g\in U^+$
si
$0=\sum (v_i\smallsmile f_i^p)(
hw)=\sum v_i(hw)f_i(hw)^p$
pour tout $h\in U^+$.

Maintenant
$\St=\Sch_{\bbF_p}(G, (p-1)\rho)^B\simeq
\Sch_{\bbF_p}(G_1B, (p-1)\rho)^B\simeq
\Sch_{\bbF_p}(U_1^+, (p-1)\rho)$
et
$\nabla(m\lambda)=\Sch_{\bbF_p}(G, m\lambda)^B\leq
\Sch_{\bbF_p}(U^+, m\lambda)$
via les restrictions.
On peut alors regarder 
$\St=\coprod_{\nu\in[0,p[^{R^+}}\bbF_p
t^\nu$
et
$\{f^p\mid f\in\nabla(m\lambda)\}
\subseteq
\coprod_{\nu\in\bbN^{R^+}}\bbF_p
t^{p\nu}$
dans l'alg\`ebre de polyn\^omes 
$\bbF_p[t]$
en les ind\'etermin\'ees  $t_{\alpha}$,
$\alpha\in R^+$,
dans laquelle on \'ecrit
$t^\nu$ pour
$\prod_{\alpha\in R^+}t_{\alpha}^{\nu_\alpha}$.
Supposons que
$\sum v_if_i^p=0$ sur  $U^+w$.
Comme $v_i(U^+w)=v_i(ww^{-1}U^+w)=(w^{-1}v_i)(\underset{\substack{\alpha\in R^+
\\
w\alpha>0}}{\prod}U_{\alpha})$
et de m\^eme pour $f_i$,
$\sum(w^{-1}v_i)(w^{-1}f_i)^p=0$
sur
$\underset{\substack{\alpha\in R^+
\\
w\alpha>0}}{\prod}U_{\alpha}$.
On peut choisir les  $v_i$ \`a partir d'une base de
$\St$ incluant
$v_-$
et tels que tous les   $v_i$ autres que 
$v_-$ soient annul\'es par  $v_-^*$.
Supposons maintenant que
$w^{-1}v_-=0$ sur $\underset{\substack{\alpha\in R^+
\\
w\alpha>0}}{\prod}U_{\alpha}$.
Alors, vu comme \'el\'ement de
$\coprod_{\nu\in[0,p[^{R^+}}\bbF_p
t^\nu$,
$w^{-1}v_-\in(t_\alpha\mid
\alpha\in R^+, w\alpha<0)$.
Mais alors $w^{-1}v_-$ aurait pour poids   
$(p-1)\rho-\underset{\substack{\alpha\in R^+\\ w\alpha<0}}{\sum}r_\alpha\alpha-\underset{\substack{\alpha\in R^+\\ w\alpha>0}}{\sum}r'_\alpha\alpha$
pour   certains 
$r_\alpha, r'_\alpha\in[0,p[$ sans qu'on ait tous les  $r_\alpha=0$.
Il s'ensuivrait que 
\[
-w^{-1}(p-1)\rho=(p-1)\rho-\underset{\substack{\alpha\in R^+\\ w\alpha<0}}{\sum}r_\alpha\alpha-\underset{\substack{\alpha\in R^+\\ w\alpha>0}}{\sum}r'_\alpha\alpha,
\]
et donc
\[
\underset{\substack{\alpha\in R^+\\ w\alpha<0}}{\sum}r_\alpha\alpha
+\underset{\substack{\alpha\in R^+\\ w\alpha>0}}{\sum}r'_\alpha\alpha
=
(p-1)(\rho+w^{-1}\rho)
=
(p-1)\underset{\substack{\alpha\in R^+\\ w\alpha>0}}{\sum}\alpha.
\]
Alors on aurait
$0>w\underset{\substack{\alpha\in R^+\\ w\alpha<0}}{\sum}r_\alpha\alpha=
w\underset{\substack{\alpha\in R^+\\ w\alpha>0}}{\sum}(p-1-r'_\alpha)\alpha\geq0$,
ce qui est absurde.
Donc, si $f_-$ est associ\'e \`a  $v_-$  dans la somme
$\sum(v_i\smallsmile f_i^p)$,
on doit avoir, comme
$\bbF_p[U^+]$ est int\`egre 
et comme les   
$w^{-1}v_i\in\coprod_{\nu\in[0,p[^{R^+}}\bbF_p
t^\nu$
sont lin\'eairement ind\'ependants des  $w^{-1}f_i^p\in\coprod_{\nu\in\bbN^{R^+}}\bbF_p
t^{p\nu}$,
 $w^{-1}f_-^p=0$ sur $\underset{\substack{\alpha\in R^+
\\
w\alpha>0}}{\prod}U_{\alpha}$,
et donc aussi  la m\^eme chose pour $w^{-1}f_-$.
Ainsi
$\sum
v_-^*(v_i)f_i(gw)=f_-(gw)=0$
pour tout $g\in U^+$,
comme voulu.
De m\^eme, on peut utiliser le scindage induit par  un vecteur de plus haut poids  $v_+ \in \St$ pour scinder de mani\`ere compatible tous les  
$X(w)$ (mais, dans l'un et l'autre cas, pas les   $X(w)$ et
$X^+(w)$ en m\^eme temps).

Finalement, montrons directement que 
 $\sigma\circ F_*v_-$
scinde de mani\`ere compatible tous les  $X(w)$, $w\in W$.
Cela revient \`a v\'erifier que si pour   tout $\tilde f\in(\Dist(U)wv_-)^\perp$ on \'ecrit 
 $v_-\smallsmile f=\sum v_i\smallsmile f_i^p$
dans  $\nabla((p-1)\rho+pm\lambda)$
comme ci-dessus, alors $\tilde f_-\in(\Dist(U)wv_-)^\perp$.
Comme les  $v_i$ forment une base de $\St$ et comme
$v_i\in\sum_{\nu\in[0,p[^{R^+}}\bbF_pt^\nu$,
les $f_i\in\sum_{\nu\in\bbN^{R^+}}\bbF_pt^\nu$
sont d\'etermin\'es de mani\`ere unique.
Supposons juste que
$\tilde f_-\ne0$ sur
$Uwv_-$.
Alors $w^{-1}f_-\not\in(t_\alpha\mid
\alpha\in R^+, w\alpha>0)$.
Mais alors $w^{-1}f\not\in(t_\alpha\mid
\alpha\in R^+, w\alpha>0)$, contredisant le fait que
$\tilde f\in(\Dist(U)wv_-)^\perp$.
 
Terminons par une question  : il est d\'emontr\'e dans 
\cite{K95} and \cite{J}  que $\cL(-\rho)$ est un facteur direct de $F_*\cO_\cB$. N\'eanmmoins, une analyse plus pr\'ecise de 
$F_*\cO_\cB$ fait apparaitre $L((p-2)\rho) \otimes \cL(-\rho)$ comme quotient de $F_*\cO_\cB$. Il est alors naturel de se demander 
si le formalisme d\'evelopp\'e ci-dessus permet de voir si  $L((p-2)\rho) \otimes \cL(-\rho)$ est facteur direct  de $F_*\cO_\cB$.

\end{document}